\pdfoutput=1
\documentclass[reqno]{amsart}

\include{EigHomotopyV9Submit.aux}
\usepackage{graphicx}
\usepackage{amssymb}
\usepackage{amsmath}
\usepackage{amsthm}
\usepackage{hyperref}

\newcommand{\abs}[1]{\left|#1\right|}

\newcommand{\im}[1]{\mbox{im}#1}
\newcommand{\Z}{\mathbb{Z}}

\newcommand{\R}{\mathbb{R}}

\newcommand{\lbreak}{\\}
\newcommand{\pgbreak}{\pagebreak}
\newcommand{\clpg}{\clearpage}
\long\def\symbolfootnote[#1]#2{\begingroup%
\def\thefootnote{\fnsymbol{footnote}}\footnote[#1]{#2}\endgroup}

\begin{document}

\footnotetext[1]{Research supported by the National Science Foundation through the Research Experiences for Undergraduates Program at Cornell.}
\footnotetext[2]{Research supported in part by the National Science Foundation, grant DMS-0652440.}

\title{Homotopies of Eigenfunctions and the Spectrum of the Laplacian on the Sierpinski Carpet}
\author{Steven M. Heilman$^{1}$}

\author{Robert S. Strichartz$^{2}$}
\thanks{Department of Mathematics, Cornell University, Ithaca, NY 14850-4201
\\
\textit{E-mail address}: \texttt{smh82@cornell.edu}
\\
Department of Mathematics, Cornell University, Ithaca, NY 14850-4201
\\
\textit{E-mail address}: \texttt{str@math.cornell.edu}
}

\date{\today}

\begin{abstract}
Consider a family of bounded domains $\Omega_{t}$ in the plane (or more generally any Euclidean space) that depend analytically on the parameter $t$, and consider the ordinary Neumann Laplacian $\Delta_{t}$ on each of them.  Then we can organize all the eigenfunctions into continuous families $u_{t}^{(j)}$ with eigenvalues $\lambda_{t}^{(j)}$ also varying continuously with $t$, although the relative sizes of the eigenvalues will change with $t$ at crossings where $\lambda_{t}^{(j)}=\lambda_{t}^{(k)}$.  We call these families \textit{homotopies} of eigenfunctions.  We study two explicit examples.  The first example has $\Omega_{0}$ equal to a square and $\Omega_{1}$ equal to a circle; in both cases the eigenfunctions are known explicitly, so our homotopies connect these two explicit families.  In the second example we approximate the Sierpinski carpet starting with a square, and we continuously delete subsquares of varying sizes.  (Data available in full at \verb!www.math.cornell.edu/~smh82!)
\end{abstract}

\maketitle


\section{Introduction}
\label{secintro}

Consider a family of bounded domains $\Omega_{t}$ in the plane (or more generally Euclidean space or a Riemannian manifold) depending on a real parameter $t$.  For each domain we consider the standard Laplacian $\Delta$ equipped with Neumann boundary conditions.  By a \textit{homotopy of eigenfunctions} we mean a continuous family of eigenfunctions $u_{t}(x)$ with eigenvalues $\lambda(t)$ (also varying continuously in $t\in\R$.
\begin{equation}\label{one1}
\begin{split}
-\Delta u_{t} & =\lambda(t)u_{t}\mbox{ on }\Omega_{t}\\
\partial_{n}u_{t} & =0\mbox{ on }\partial\Omega_{t}
\end{split}
\end{equation}
For a reasonable domain $\Omega_{t}$ there is a complete orthonormal basis of eigenfunctions.  We do not expect every eigenfunction to belong to a homotopy.  If the eigenvalue is simple, then it is easy to see that the (unique up to a constant) eigenfunction does belong to a homotopy for at least a small interval in the parameter $t$.  However, if the eigenspace has higher multiplicity at a fixed value of $t$, the best we can expect is that there is a special basis of eigenfunctions that belong to a homotopy.

We face the same situation if we consider a continuous family $A_{t}$ of symmetric $n\times n$ matrices.  The eigenvalues depend continuously on $t$, but the eigenvectors do not necessarily vary continuously.  For example, let
\begin{equation}\label{one2}
A_{t}=
\begin{cases}
t^{3}B_{1} & t\geq0\\
t^{3}B_{2} & t<0
\end{cases}
\end{equation}
for $B_{1},B_{2}$ $2\times 2$ symmetric matrices with different eigenvectors.  One can avoid this pathology by assuming that the derivative with respect to $t$ is nonzero, or that the dependence on $t$ is analytic.  However, if $A_{t}=tB$ where $B$ has distinct eigenvectors, then for $t=0$, the matrix has eigenspace corresponding to $\lambda=0$ containing all vectors, whereas only the eigenvectors of $B$ belong to homotopies.

Observe now that the double eigenvalue of $A_{t}$ at $t=0$ is unstable under perturbation.  More specifically, recall that matrices with multiple eigenvalues, considered as a subspace of $n\times n$ symmetric matrices, have codimension $2$.  Similarly, those matrices with multiple eigenvalues in $n\times n$ Hermitian matrices have codimension $3$.  So, for a generic one-parameter family of self-adjoint matrices, multiple eigenvalues are very rare.  Moreover, in the self-adjoint case, if two eigenvalues move near each other (say, as $t$ increases), the eigenvalues then repel each other.  We can see this from the following equation (derived below as in the article ``When are eigenvalues stable?'' \cite{tao09}):
\begin{equation}\label{one2.5}
\lambda''_{k}=u_{k}^{*}A_{t}''u_{k}
+2\sum_{j\neq k}\frac{\abs{u_{k}^{*}A_{t}'u_{j}}}{\lambda_{k}-\lambda_{j}}
\end{equation}
where $u_{k}$ is the $k^{th}$ eigenfunction (an $L^{2}$ normalized column vector), $\lambda_{j}$ is the $j^{th}$ eigenvalue, $^{*}$ denotes complex transpose, and $'$ denotes a derivative in the $t$ variable.  Note that the eigenvectors and eigenvalues are functions of $t$.

To obtain Eq. (\ref{one2.5}), we differentiate the equation $A_{t}u_{k}(t)=\lambda_{k}(t)u_{k}(t)$ in the $t$ variable to get
\begin{equation}\label{one2.55}
A_{t}u_{k}'+A_{t}'u_{k}=\lambda_{k}u_{k}'+\lambda_{k}'u_{k}
\end{equation}
We then differentiate the identity $u_{k}^{*}u_{k}=1$ to get $u_{k}^{*}u_{k}'=0$.  Then, applying $u_{k}^{*}$ to both sides of Eq. (\ref{one2.55}) and using the previous sentence gives
\begin{equation}\label{one2.6}
\lambda_{k}'=u_{k}^{*}A_{t}'u_{k}
\end{equation}
Notice we have also used the fact that
\begin{equation}\label{one2.61}
u_{k}^{*}A_{t}=\lambda_{k}u_{k}^{*},
\end{equation}
which is just the eigenvalue equation for $u_{k}$ with the adjoint taken on both sides.  Now, apply $u_{j}^{*}$ to both sides of Eq. (\ref{one2.55}) for $j\neq k$, and apply Eq. (\ref{one2.61}) and the orthogonality of the eigenvectors to the result to obtain
\begin{equation}\label{one2.63}
u_{j}^{*}A_{t}'u_{k}+(\lambda_{j}-\lambda_{k})u_{j}^{*}u_{k}'=0.
\end{equation}
So, re-writing, summing over $j$, using the completeness of the eigenvectors (and recalling that $u_{j}^{*}u_{j}'=0$) gives
\begin{equation}\label{one2.65}
u_{k}'=\sum_{j\neq k}\frac{u_{j}^{*}A_{t}'u_{k}}{\lambda_{k}-\lambda_{j}}
\end{equation}
Now, differentiating Eq. (\ref{one2.6}) gives
\begin{equation}\label{one2.67}
\lambda_{k}''=u_{k}^{*}A_{t}''u_{k}+u_{k}^{*}A_{t}'u_{k}'+(u_{k}^{*})'A_{t}'u_{k}
\end{equation}
And after applying Eq. (\ref{one2.65}) to Eq. (\ref{one2.67}), we finally get Eq. (\ref{one2.5}).

We now wish to comment further on the significance of Eq. (\ref{one2.5}).  The repulsion of two distinct eigenvalues, or level repulsion, manifests in the second term on the right of Eq. (\ref{one2.5}).  That is, if we increase $t$ linearly, and two eigenvalues $\lambda_{k}$ and $\lambda_{k'}$ move close together, the second term of Eq. (\ref{one2.5}) becomes very large (assuming the numerator does not vanish).  In fact, $\lambda_{k}''$ has a term of the form $1/(\lambda_{k}-\lambda_{k'})$, while $\lambda_{k'}''$ has a term of the form $1/(\lambda_{k'}-\lambda_{k})$.  So, the larger of the two eigenvalues is forced upward, and the smaller of the two is forced downward.  Therefore, Eq. (\ref{one2.5}) accounts for generic eigenvalue repulsion, in this finite dimensional case.

Now, can these observations carry over to the infinite dimensional setting?  We are concerned here with eigenvalue repulsion and rarity of multiple eigenvalues.  Specifically, do these effects occur for compact planar domains which vary with a time parameter?  To a certain extent, yes.  Restricting an infinite dimensional operator to finitely many dimensions as in Eq. (\ref{one5}) below, the repulsion effects from the finite dimensional Hermitian case should carry over to self-adjoint operators.  (For rigorous results on the rarity of multiple eigenvalues, see \cite{verdiere88,lamberti06,teytel99}.)

Mimicking the finite dimensional case, we reformulate our original problem by taking a fixed domain $\Omega$ and varying the differential operator.  We assume that we have a family of diffeomorphisms $\Psi_{t}:\Omega\to\Omega_{t}$ that depends analytically on the parameter $t$.  Then we can pull back the Euclidean metric on $\Omega_{t}$ to a Riemannian metric on $\Omega$.  The pullback of the standard Laplacian on $\Omega_{t}$ becomes the Laplace-Beltrami operator on $\Omega$, which we denote by $\Delta_{t}$.  Similarly, we can pull back the energy and measure on $\Omega_{t}$ to the corresponding quantities on $\Omega$ for the Riemannian metric.  We note that the Neumann boundary conditions correspond, and also the coefficients of $\Delta_{t}$ are analytic functions in $t$.  Let $[a,b]$ be any finite interval in $[0,\infty)$ such that the endpoints are not Neumann eigenvalues of any $\Delta_{t}$ (to do this we may have to restrict the interval in which $t$ varies).  Then the spectral projection operator $P_{t}([a,b])$ associated with $-\Delta_{t}$ and the interval $[a,b]$ is well-defined and of finite rank, and may be expressed as a contour integral in the resolvent
\begin{equation}\label{one3}
R(\lambda)=\int_{\Gamma}(\lambda+\Delta_{t})^{-1}d\lambda
\end{equation}
for a suitable contour $\Gamma$ in the complex plane enclosing the interval.  The resolvent depends analytically on $t$ (this may be seen via the pseudodifferential calculus, for example), so the projection operator also depends analytically \textbf{on} $t$.  If the rank of the projection operator is $N$, we can find an $N$-dimensional space $V$ of functions such that $P_{t}([a,b])$ maps $V$ onto its range.  Let $P_{t}^{-1}([a,b])$ denote the inverse.  Then we can transform the eigenvalue equation
\begin{equation}\label{one4}
-\Delta_{t}u_{t}=\lambda(t)u_{t}\mbox{ on }\Omega
\end{equation}
into the eigenvalue equation
\begin{equation}\label{one5}
-\left(P_{t}([a,b])^{-1}\Delta_{t}P_{t}([a,b])\right)x_{t}=\lambda(t)x_{t}\mbox{ in }V
\end{equation}
via $x_{t}=P_{t}([a,b])^{-1}u_{t}$.  Now we are in the finite dimensional setting with a family of operators depending analytically on the parameter, so $N$ homotopies of eigenfunctions exist, and these yield $N$ homotopies of eigenfunctions $u_{t}$ on $\Omega$.

Once we have the existence of homotopies for small intervals of the parameter and bounded regions of the spectrum, we can combine them to obtain global homotopies for the whole spectrum.  If $u_{t_{0}}$ is an eigenfunction corresponding to a simple eigenvalue $\lambda(t_{0})$ of $\Delta_{t_{0}}$, then there is a unique homotopy $u_{t}$ with $u_{t}|_{t=t_{0}}=u_{t_{0}}$.  If $\lambda(t_{0})$ is an eigenvalue of $\Delta_{t_{0}}$ with nontrivial multiplicity, then its eigenspace $E(t_{0})$ has a basis of eigenvectors which lie in homotopies.  In fact one can take this basis to be orthonormal.  Suppose for example, that $\lambda(t_{0})$ has multiplicity $2$, and let $u_{t}^{(1)}$ and $u_{t}^{(2)}$ denote $2$ homotopies with eigenvalues $\lambda^{(1)}(t)$ and $\lambda^{(2)}(t)$.  Then either $\lambda^{(1)}(t)=\lambda^{(2)}(t)$ for all $t$ or $\lambda^{(1)}(t)\neq\lambda^{(2)}(t)$ for $0<\abs{t-t_{0}}<\epsilon$.  In the first case we may take linear combinations of $u_{t}^{(1)}$ and $u_{t}^{(2)}$ to obtain orthogonality.  In the second case we have the orthogonality of $u_{t}^{(1)}$ and $u_{t}^{(2)}$ for $0<\abs{t-t_{0}}<\epsilon$ (because the eigenvalues are distinct), and by continuity it must hold for $t=t_{0}$.

The deeper question, however, is what is the significance of the homotopies?  If we have homotopies between domains $\Omega_{0}$ and $\Omega_{1}$ then we have a pairing up of the eigenfunctions on $\Omega_{0}$ and $\Omega_{1}$.  What does this tell us about the paired eigenfunctions?  The goal of this paper is to investigate this question experimentally.  We examine two different examples.  In the first, something of a warmup exercise, we take $\Omega_{0}$ be a circle and $\Omega_{1}$ to be a square.  We know exactly what the eigenfunctions are at both endpoints.  For the unit circle, using polar coordinates, the eigenspaces typically have multiplicity $2$ and are spanned by $J_{m}(\lambda_{mk}r)\cos m\theta$ and $J_{m}(\lambda_{mk}r)\sin m\theta$ where the eigenvalues $\lambda_{mk}$ are the successive zeroes of $J_{m}'$, indexed by $k=1,2,3,\ldots$.  Here $m=1,2,3,\ldots$.  The exceptional case $m=0$ yields a $1$-dimensional eigenspace spanned by $J_{0}(\lambda_{0k}r)$.  For the unit square, there is a basis of eigenfunctions of the form $\cos\pi jx\cos\pi ky$ for $j=0,1,2\ldots$ and $k=0,1,2,\ldots$ with eigenvalues $\pi^{2}(j^{2}+k^{2})$.  Typical eigenspaces have dimension $1$ ($j=k$) or $2$ ($j\neq k$), but higher multiplicity may occur due to coincidences, such as $3^{2}+4^{2}=5^{2}+0^{2}$ or $7^{2}+1^{2}=5^{2}+5^{2}$.

In order to simplify and clarify the discussion we consider the action of the dihedral symmetry group $D_{4}$ on the eigenspaces.  $D_{4}$ acts on the circle and square and all the domains we consider.  It is well known that this $8$ element group has $5$ irreducible representations, one being $2$-dimensional and the others $1$-dimensional.  Each eigenspace either transforms according to one of these representations or splits into a direct sum of such pieces, and the homotopies respect the representation types.  We denote the $1$-dimensional representations $1++,1+-,1-+,$ and $1--$, where the first $\pm$ refers to the symmetry or skew-symmetry with respect to diagonal reflections, and the second $\pm$ refers to the symmetry or skew-symmetry with respect to horizontal and vertical reflections.  We now describe the eigenspaces sorted according to representation type for the circle and the square.

\underline{Circle:}

$1++:$ The eigenspace $E_{++}^{c}(4m,k)$ is generated by\lbreak $J_{4m}(\lambda_{(4m)k}r)\cos 4m\theta$ for $m=0,1,2,\ldots$ and $k=1,2,3,\ldots$

$1-+:$ The eigenspace $E_{-+}^{c}(4m,k)$ is generated by\lbreak $J_{4m}(\lambda_{(4m)k}r)\sin 4m\theta$ for $m=1,2,3,\ldots$ and $k=1,2,3,\ldots$

$1+-:$ The eigenspace $E_{+-}^{c}(4m+2,k)$ is generated by\lbreak $J_{4m+2}(\lambda_{(4m+2)k}r)\cos (4m+2)\theta$ for $m=0,1,2,\ldots$ and $k=1,2,3,\ldots$

$1--:$ The eigenspace $E_{--}^{c}(4m+2,k)$ is generated by\lbreak $J_{4m+2}(\lambda_{(4m+2)k}r)\sin (4m+2)\theta$ for $m=0,1,2,\ldots$ and $k=1,2,3,\ldots$

$2:$ The eigenspace $E_{2}^{c}(2m+1,k)$ is generated by\lbreak $J_{2m+1}(\lambda_{(2m+1)k}r)\cos (2m+1)\theta$ and $J_{2m+1}(\lambda_{(2m+1)k}r)\sin (2m+1)\theta$ for $m=0,1,2,\ldots$ and $k=1,2,3,\ldots$

\underline{Square:}

$1++:$ The eigenspace $E_{++}^{s}(2j,2k)$ is generated by\lbreak $\cos\pi2jx\cos\pi2ky+\cos\pi2kx\cos\pi2jy$ for $0\leq j\leq k$

$1-+:$ The eigenspace $E_{-+}^{s}(2j,2k)$ is generated by\lbreak $\cos\pi2jx\cos\pi2ky-\cos\pi2kx\cos\pi2jy$ for $0\leq j< k$

$1+-:$ The eigenspace $E_{+-}^{s}(2j+1,2k+1)$ is generated by\lbreak $\cos\pi(2j+1)x\cos\pi(2k+1)y+\cos\pi(2k+1)x\cos\pi(2j+1)y$ for $0\leq j\leq k$

$1--:$ The eigenspace $E_{--}^{s}(2j+1,2k+1)$ is generated by\lbreak $\cos\pi(2j+1)x\cos\pi(2k+1)y-\cos\pi(2k+1)x\cos\pi(2j+1)y$ for $0\leq j< k$

$2:$ The $2$-dimensional eigenspace $E_{2}^{s}(2j+1,2k)$ is generated by\lbreak $\cos\pi(2j+1)x\cos\pi2ky$ and $\cos\pi2kx\cos\pi(2j+1)y$ for $0\leq j\leq k$

In the case of the square, because of coincidences, the full eigenspace associated with a given eigenvalue may be a direct sum of some of the spaces described above, for example $E_{2}^{s}(5,0)\oplus E_{2}^{s}(3,4)$ or $E_{+-}^{s}(5,5)\oplus E_{+-}^{s}(7,1)$.  When we construct homotopies, we often find that they do not pass through the individual summands.  This problem does not arise for the circle, since as far as we know there are no coincidences among the eigenvalues $\lambda_{mk}$.

In section \ref{secctos} we describe the homotopies of eigenfunctions for two different analytic families of domains joining the circle to a square.  In one, the circle expands out to a square in horizontal/vertical orientation, and in the other the circle contracts in to a square in diagonal orientation.  (The fact that these are not unit squares just requires a trivial adjustment.)  In both cases the homotopies connect the same $E^{s}$ and $E^{c}$ eigenspaces.  We observe a canonical splitting of two-dimensional eigenspaces at the endpoints of the homotopy.  We also do not find any two-dimensional eigenspaces away from the endpoints of the homotopy.

In our second example in section \ref{seccarp} we consider a sequence of domains $A_{i}$ that approximate the Sierpinski carpet ($SC$).  Let $F_{1},\ldots, F_{8}$ denote the similarities that map to square $A_{0}$ into the squares of the ``picture frame'' $A_{1}$ that is $A_{0}$ with the middle ninth removed.  Then
\begin{equation}\label{one6}
SC=\cup_{i=1}^{8}F_{i}(SC)
\end{equation}
and
\begin{equation}\label{one7}
A_{j}=\cup_{i=1}^{8}F_{i}(A_{j-1})
\end{equation}
In a previous work \cite{berry08} we gave experimental evidence that the spectra of the Neumann Laplacians on $A_{j}$ converge, after the appropriate renormalization, to the spectrum of the Neumann Laplacian on $SC$, with the eigenfunctions also converging.  In other words, there is a complete orthonormal basis $u_{n}^{(j)}$ of eigenfunctions
\begin{equation}\label{one8}
-\Delta u_{n}^{(j)}=\lambda_{n}^{(j)}u_{n}^{(j)}
\end{equation}
on $A_{j}$ with Neumann boundary conditions such that
\begin{equation}\label{one9}
\lim_{j\to\infty}r^{j}\lambda_{n}^{(j)}=\lambda_{n}\mbox{ and}
\end{equation}
\begin{equation}\label{one10}
\lim_{j\to\infty}u_{n}^{(j)}|_{SC}=u_{n}
\end{equation}
for some $r\in\R$.  We would then identify $u_{n}$ as the Neumann eigenfunctions of the Barlow-Bass Laplacian $\Delta$ on $SC$ \cite{barlow89}.  The recently announced proof \cite{barlow08} of the uniqueness of that Laplacian justifies this identification.  Combining this uniqueness result, the resolvent convergence \cite{bass97} (see \S7), and the relation $e^{t\Delta}f(x)=E_{x}[f(X_{t})]$ between the Laplacian $\Delta$ and Brownian motion $X_{t}$, suggest a method of rigorous proof of Eq. (\ref{one10}).  The relation $e^{t\Delta}f(x)=E_{x}[f(X_{t})]$ is important because it shows that the speeding up of the Brownian motion \cite{bass97} exactly corresponds to multiplication of $\Delta$ by a constant.

However, even with convergence in Eq. (\ref{one10}), there is an obstacle to uncovering the dynamics of the eigenvalues and eigenfunctions.  It was noted in \cite{berry08} that, for a given $n$, the eigenvalues $\{\lambda_{n}^{(j)}\}$ are not necessarily in increasing order (with respect to the original ordering for $j=0$).  In order to better understand the discrete sequence $u_{n}^{(j)}$ of eigenfunctions, we embed them in homotopies.  For example, to connect the square $A_{0}$ with the picture frame $A_{1}$, we first remove the center of the square.  Since points have zero capacity in the square, this does not change the Neumann Laplacian.  We then gradually enlarge the hole by deleting squares in the center with side length going from $0$ to $1/3$.  Similarly, to go from $A_{1}$ to $A_{2}$ we enlarge $8$ holes in the center of the squares $F_{i}A_{0}$ with side length going from $0$ to $1/9$.  We carried out the computation of these homotopies up to $A_{4}$.  (Note that it would also be possible to do a single deformation from $A_{0}$ to $A_{4}$ by opening all holes simultaneously.)  Some eigenfunctions change only slightly during the course of the homotopy; we could say that they make minor adjustments to the changing domains.  However, most eigenfunctions show considerable variation, especially for higher eigenvalues.  Ideally we would hope to make a reasonable guess from the shape of $u_{n}^{(4)}$ as to which $u_{n'}^{(0)}$ it is connected to.  We are very far from being able to do this.  In fact, the following cautionary tale shows some of the obstacles.  If $j$ and $k$ are divisible by $3$, then all the eigenfunctions in any of the $E^{s}(j,k)$ spaces are also Neumann eigenfunctions on $A_{1}$ (when restricted to $A_{1}$).  However, many of the homotopies do not connect the restriction to $A_{1}$ with the original eigenfunction on $A_{0}$.

\section{Circle to Square Homotopy}
\label{secctos}

Let $L\subset\R^{2}$ denote the line $x+y=1$ in the closed first quadrant.  In polar coordinates $L$ is parameterized by the equation $r=1/(\cos\theta+\sin\theta)$, $\theta\in[0,\pi/2]$.  Then, the map $\psi(r,\theta)=(1/(\cos\theta+\sin\theta)r,\theta)$ defined in the closed first quadrant takes the arc of the circle of radius $1$ centered at the origin to the line $L$.  In fact, it also maps the arc of the circle of radius $r$ to the line $\{(x,y):x+y=r\}$.  With this in mind, let $H(t,r,\theta):[0,1]\times\R^{2}\to\R^{2}$ be the homotopy
\begin{equation}\label{two1}
H(t,r,\theta)=\left(\left((1-t)+\frac{t}{\cos\theta+\sin\theta}\right)r,\theta\right)
\end{equation}
For $\theta\in[0,\pi/2]$, this is the straight-line homotopy between the identity map of $\R^{2}$ and $\psi(r,\theta)$.  So, as we vary $t$, $H$ restricted to the closed first quadrant maps the square diamond of side length $\sqrt{2}$ to the unit disc.  Observe that for $t\in(0,1)$, $H$ (with the same restriction) maps the square to a half-ellipse.

Rotating the domain and range of the map into all four quadrants, we get a straight-line homotopy which connects the square diamond to the unit disc and which is analytic in the time parameter $t$ on $(0,1)$ (and continuous on $[0,1]$).  If we view $H(t,r,\theta)$ as an analytic deformation of the standard metric on the unit disc, the Neumann Laplacian $\Delta_{t}$ on the disc is a holomorphic family of type (A) (see \cite{kato66}).  Then, Theorem 3.9,VII.\S4.5 from the same book says the eigenvalues and eigenvectors can be represented as functions holomorphic in $t$.  More specifically, the eigenvalues are holomorphic functions and the eigenfunctions can be chosen to be holomorphic with a Banach space as their range.  In our terminology, we have an (analytic) homotopy of eigenfunctions.  This is the theoretical basis of our investigation.

An analogous process can be achieved for $\phi(r,\theta)=(\sqrt{2}/(\cos\theta+\sin\theta)r,\theta)$ and
\begin{equation}\label{two2}
F(t,r,\theta)=\left(\left((1-t)+\frac{\sqrt{2}t}{\cos\theta+\sin\theta}\right)r,\theta\right)
\end{equation}
where $\phi$ maps the unit circle in the first quadrant to the square diamond of side length $2$, and upon reflection as above, $F$ gives a straight line homotopy between the identity map on the unit circle and $\phi(r,\theta)$.  Again, we get a straight-line homotopy from the unit disc to the square, and with previously mentioned results \cite{kato66}, we have a homotopy of eigenfunctions.  Note that $F$ and $H$ are not identical.

Since all of our domains are invariant under the $D_{4}$ symmetry group, we can simplify the eigenfunction computations by reducing to a fundamental domain.  On this domain we impose appropriate boundary conditions according to the representation type.  For the $1$-dimensional representation, we restrict to the sector $0\leq\theta\leq\pi/4$.  Recall that the functions will extend evenly when reflected about $\theta=0$ in the $1++$ and $1-+$ cases, and oddly in the $1+-$ and $1--$ cases.  When reflecting about $\theta=\pi/4$, the functions will extend evenly in the $1++$ and $1+-$ cases, and oddly in the $1-+$ and $1--$ cases.  Note that performing an even extension across a ray is equivalent to imposing Neumann boundary conditions on that ray.  Similarly, the odd extension is equivalent to Dirichlet conditions.

For the $2$-dimensional representation our fundamental domain is the sector $0\leq\theta\leq\pi/2$, and we impose Neumann boundary conditions on the ray $\theta=0$ and Dirichlet conditions on the ray $\theta=\pi/2$.  Notice that our fundamental domains are simply connected.  Therefore, as in Theorem 6.4 of \cite{teytel99}, we should expect no two-dimensional eigenspaces along this homotopy (for each individual symmetry family).

In order for Matlab R2006a to understand our geometry, we must give the geometric data to the program in its desired format. To accomplish this task in Matlab's PDE Toolbox for a fixed $t$, we specify polygons and circles.  More specifically, we specify the vertices of the polygon and the center and radius of a circle, with standard coordinates in the plane.  Once we input this data, we then fix a closed interval of the real line in which we look for eigenvalues $\lambda$.  The program then uses the Arnoldi Algorithm with function \verb!pdeeig! to solve the eigenvalue equation $(-\Delta)u=\lambda u$ for $\lambda$ in the specified interval.  For this particular homotopy, this process yields sixty to seventy eigenvalues within minutes.  As our first approximation to a smooth homotopy using $H$ or $F$, we solve for eigenfunctions at times $t\in\{k/10\}_{k=0}^{10}$.  For now, Matlab's solution method is sufficient. However, other solvers based in different programs may be more effective in the future. For instance, the boundary element method appears much more efficient than the finite element method.

        In the following diagrams, we show the observed eigenvalue dynamics for each of the five symmetry families ($1++,1+-,1-+,1--,2$).  For convenience, we divide all eigenvalues by $\pi^{2}$.  As stated above, the data is collected for $t\in\{k/10\}_{k=0}^{10}$.  As $t$ varies, we determine whether or not two eigenvalues cross by inference (see below for more details).

        Observe that the eigenvalue plots exhibit eigenvalue repulsion, as in Eq. (\ref{one2.5}) in the finite dimensional setting.  When two eigenvalues draw near each other, they ``collide'' and then veer away rapidly (see Fig. \ref{figthirteen.5} (a) below).  Intuitively, when two eigenvalues move near each other, one could imagine a slightly perturbed geometry for which these eigenvalues would be the same.  Indeed, the eigenfunctions also support this intuition.  If we examine the eigenfunctions before and after the collision, we see that they swap identities.  That is, if we have eigenvalues $\lambda_{j+1}>\lambda_{j}$ which collide around $t=t_{0}$, then $u_{j+1}(t_{0}-\epsilon)$ looks like $u_{j}(t_{0}+\epsilon)$, and $u_{j}(t_{0}-\epsilon)$ looks like $u_{j+1}(t_{0}+\epsilon)$. In the intermediate region $(t_{0}-\epsilon,t_{0}+\epsilon)$, we then observe two linear combinations of our eigenfunctions $u_{j+1}(t_{0}-\epsilon)$ and $u_{j}(t_{0}-\epsilon)$.  This effect has been observed before \cite{backer97} and coined the superposition effect.  Notice that this superposition effect only occurs among eigenfunctions of the same symmetry class.

\begin{figure}[htbp!]
\begin{center}
\includegraphics[scale=1]{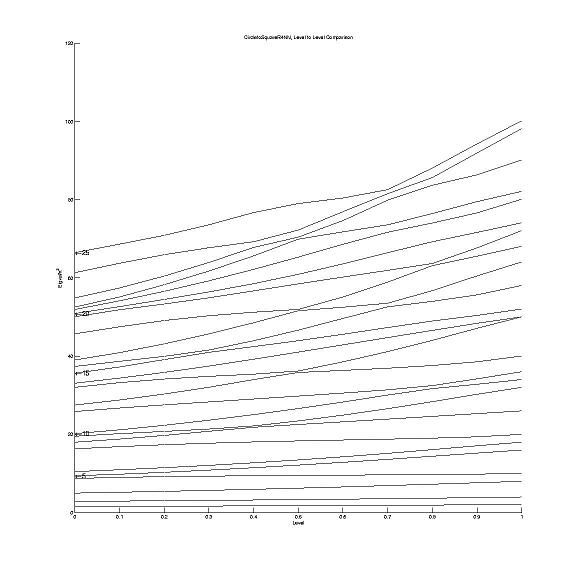}
\end{center}
\caption{$(1++)$ Eigenvalue Summary Plot.  Circle ($t=0$), Square ($t=1$)}
\label{figone}
\end{figure}

\begin{figure}[htbp!]
\begin{center}
\includegraphics[scale=1]{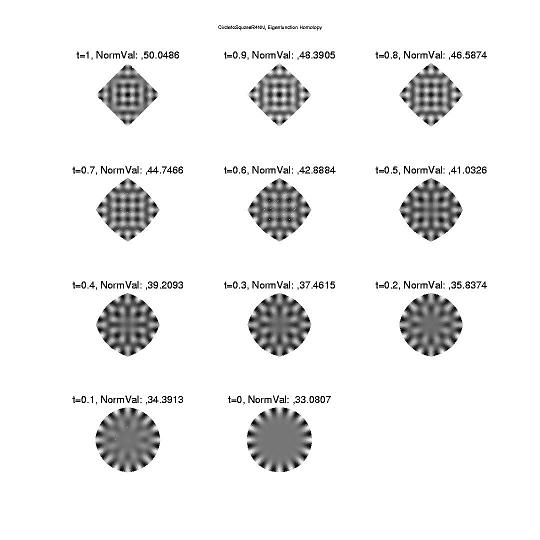}
\end{center}
\caption{$(1++)$ EigVal $14$ Homotopy}
\label{figsix.1}
\end{figure}

\begin{figure}[htbp!]
\begin{center}
\includegraphics[scale=1]{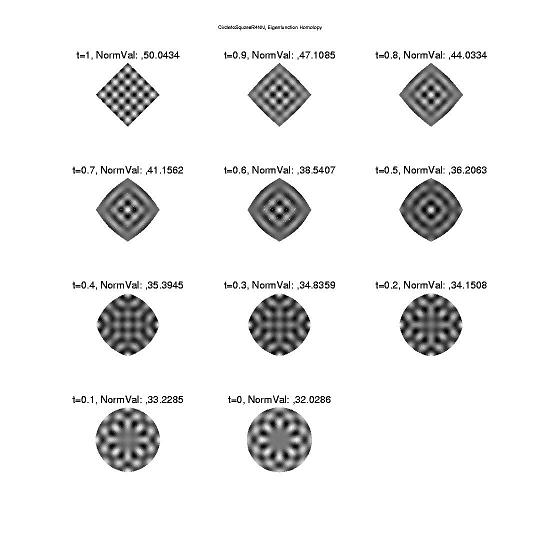}
\end{center}
\caption{$(1++)$ EigVal $13$ Homotopy}
\label{figsix}
\end{figure}

        When a two-dimensional eigenspace occurs (analytically), the superposition effect is sometimes not observed numerically.  For example, consider eigenfunctions $u_{13}$ and $u_{14}$ from the $(1++)$  family on the square.  We see from Fig. \ref{figone} that these two eigenfunctions separate from the two dimensional eigenspace at $t=1$, as $t$ decreases.  A priori, this two-dimensional eigenspace could split into any two (orthogonal) eigenfunctions from a two-parameter family $\{a_{1}u_{13}+a_{2}u_{14}\colon a_{1},a_{2}\in\R\}$.  However, we see from Figs. \ref{figsix} and \ref{figsix.1} that these eigenfunctions appear to split according to the canonical eigenfunctions on the square.  This behavior is consistently observed for this homotopy.  We will see in Section \ref{seccarp} different behavior for that homotopy.

\begin{figure}[htbp!]
\begin{center}
\includegraphics[scale=1]{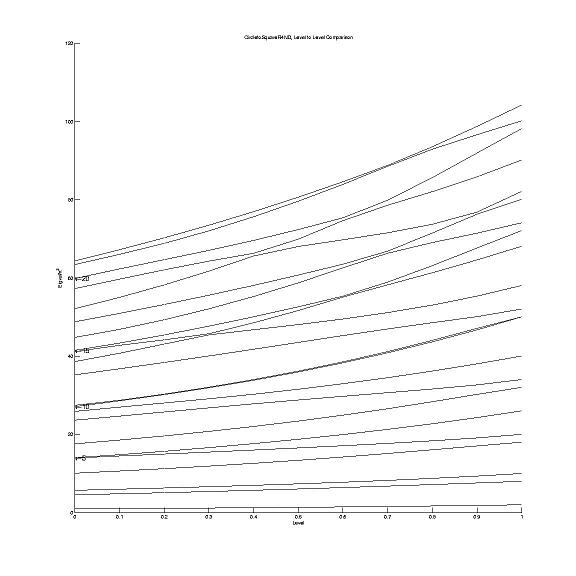}
\end{center}
\caption{$(1+-)$ Eigenvalue Summary Plot}
\label{figtwo}
\end{figure}

\begin{figure}[htbp!]
\begin{center}
\includegraphics[scale=1]{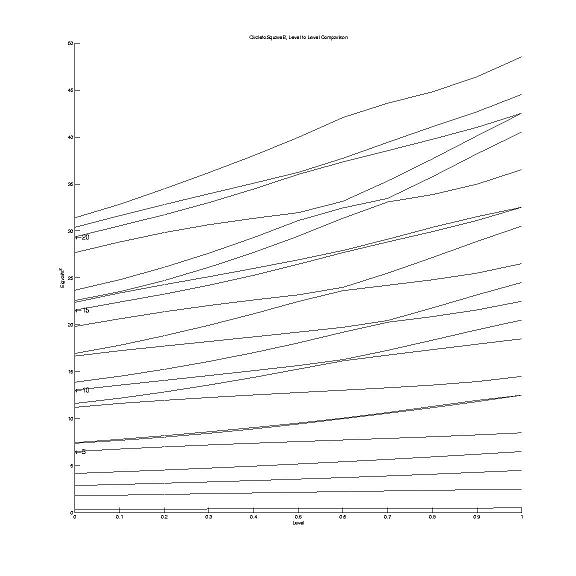}
\end{center}
\caption{$(2)$ Eigenvalue Summary Plot}
\label{figfive}
\end{figure}

        Additionally, observe the closely bound eigenfunctions in some of our eigenvalue plots.  That is, as $t$ varies in our homotopy, we see two eigenfunctions whose eigenvalues remain very close for all $t$.  This effect is puzzling for the authors to explain, and it is unclear whether it only occurs in the lower part of the spectrum.  We observe this effect for eigenfunctions $10$ and $11$ in Fig. \ref{figtwo} and eigenfunctions $5$ and $6$ in Fig. \ref{figfive}.  This effect is not observed in our homotopy of Section \ref{seccarp}.

        Also note that certain eigenvalues change in similar ways.  For instance, eigenvalues $13$,$18$, and $24$ of the circle in the $(1++)$ family all have noticeably negative slope around $t=.5$.

        This observation indicates that patterns occur in our eigenfunction homotopies.  In other words, it seems that sets of similar eigenfunctions behave in similar ways in the homotopies.

\begin{figure}[htbp!]
\begin{center}
\includegraphics[scale=1]{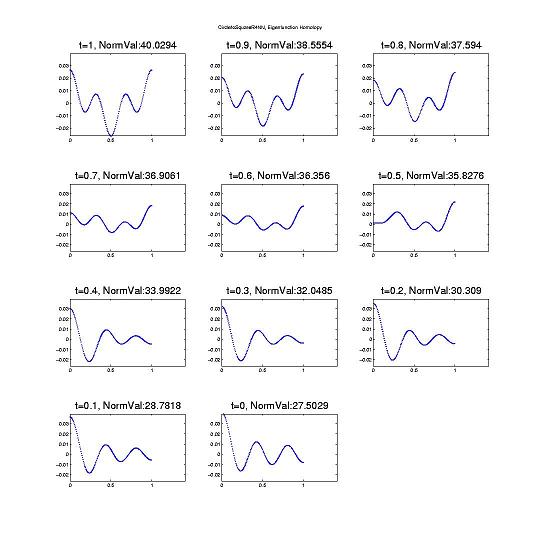}
\end{center}
\caption{$(1++)$ EigVal $13$ Homotopy, Restricted to a Line}
\label{figseven}
\end{figure}

    In Fig. \ref{figsix}, we depict the homotopy of eigenfunctions of Eigenvalue $13$ of the circle in the $(1++)$ family.  In Fig. \ref{figseven}, we show this same homotopy restricted to the positive $x$-axis.

\begin{figure}[htbp!]
\begin{center}
\includegraphics[scale=1]{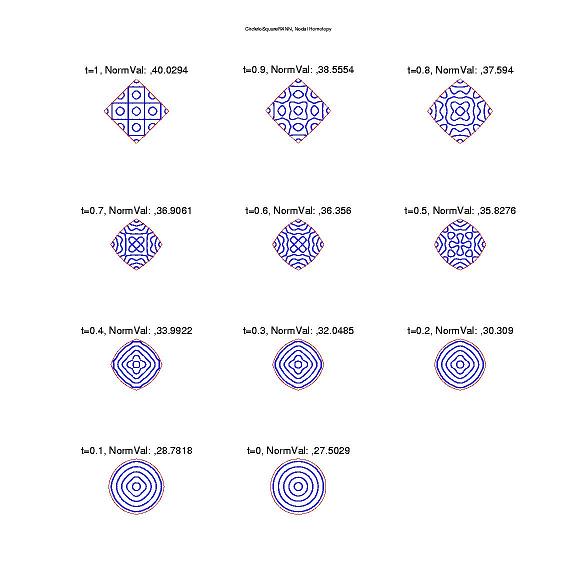}
\end{center}
\caption{$(1++)$ EigVal $13$ Nodal Homotopy}
\label{figeight}
\end{figure}

    In Fig \ref{figeight}, we take the same homotopy as in Figs. \ref{figsix} and \ref{figseven}, but for simplicity we only examine the nodal set of the eigenfunction.  We approximate this nodal set by the inverse image of a small neighborhood of $0$ of the form $(-\epsilon,\epsilon)$.

\begin{figure}[htbp!]
\begin{center}
\includegraphics[scale=1]{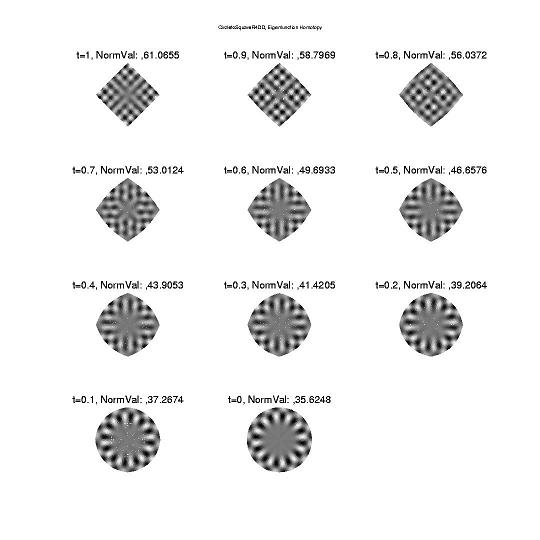}
\end{center}
\caption{$(1--)$ EigVal $10$ Homotopy}
\label{fignine}
\end{figure}

\begin{figure}[htbp!]
\begin{center}
\includegraphics[scale=1]{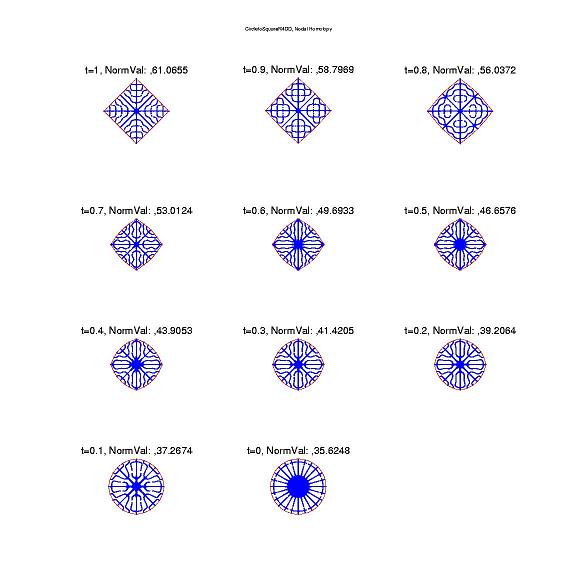}
\end{center}
\caption{$(1--)$ EigVal $10$ Nodal Homotopy}
\label{figten}
\end{figure}

    In Fig. \ref{fignine} we depict the homotopy of Eigenvalue $10$ on the circle in the $(1--)$ family.  In Fig. \ref{figten} we show the homotopy of the corresponding nodal set, as in Fig. \ref{figeight}.

\begin{figure}[htbp!]
\begin{center}
\includegraphics[scale=1]{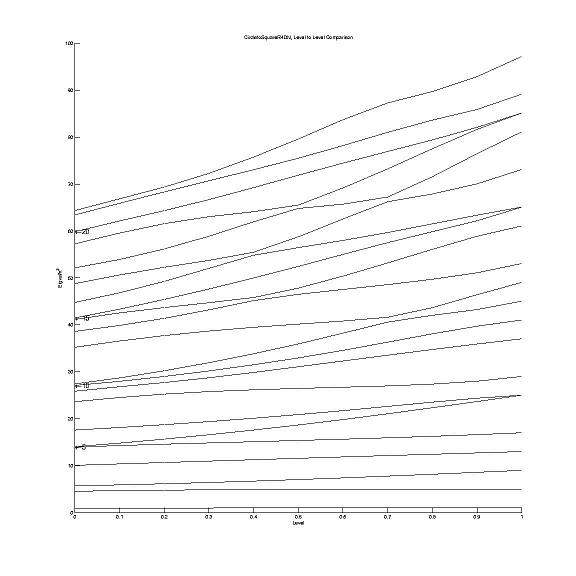}
\end{center}
\caption{$(1-+)$ Eigenvalue Summary Plot}
\label{figthree}
\end{figure}

\begin{figure}[htbp!]
\begin{center}
\includegraphics[scale=1]{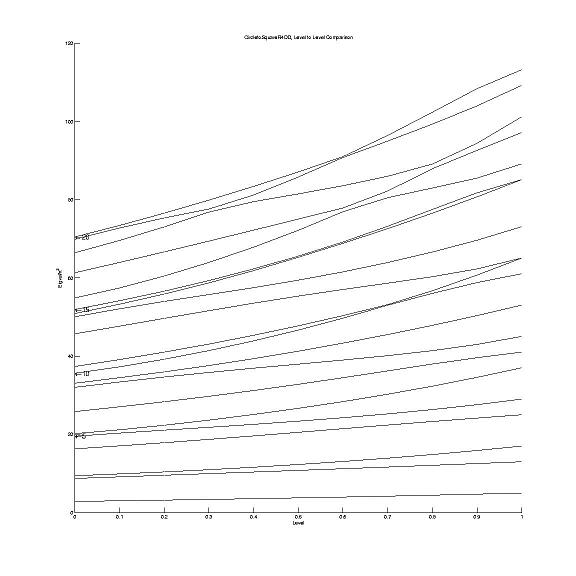}
\end{center}
\caption{$(1--)$ Eigenvalue Summary Plot}
\label{figfour}
\end{figure}

\begin{figure}[htbp!]
\begin{center}
\includegraphics[scale=1]{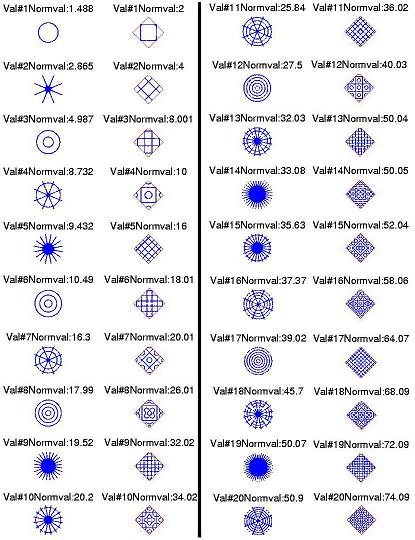}
\end{center}
\caption{$(1++)$ Eigenfunction correspondences}
\label{figeleven1}
\end{figure}

\begin{figure}[htbp!]
\begin{center}
\includegraphics[scale=1]{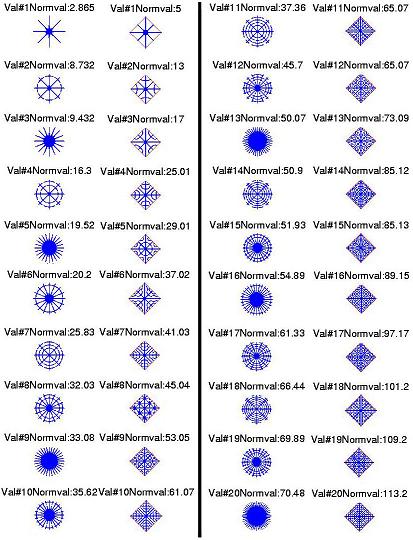}
\end{center}
\caption{$(1--)$ Eigenfunction correspondences}
\label{figeleven2}
\end{figure}

\begin{figure}[htbp!]
\begin{center}
\includegraphics[scale=1]{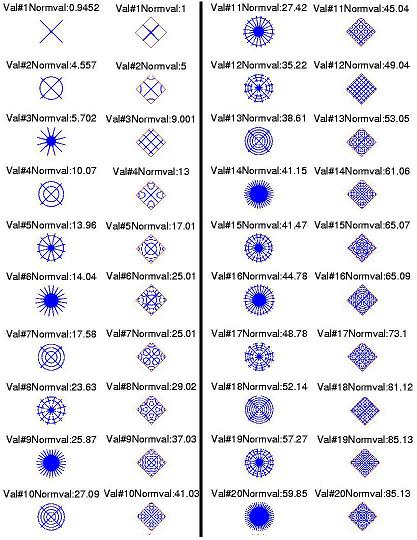}
\end{center}
\caption{$(1-+)$ Eigenfunction correspondences}
\label{figeleven3}
\end{figure}

\begin{figure}[htbp!]
\begin{center}
\includegraphics[scale=1]{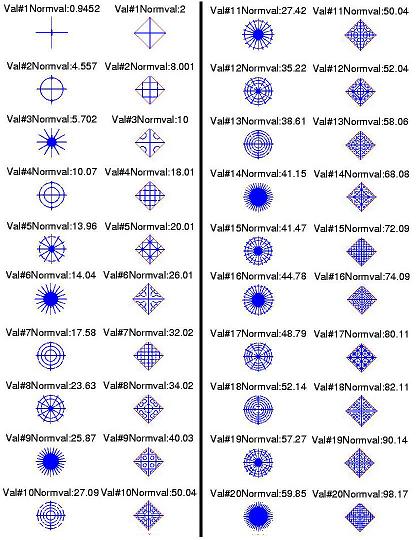}
\end{center}
\caption{$(1+-)$ Eigenfunction correspondences}
\label{figeleven4}
\end{figure}

\begin{figure}[htbp!]
\begin{center}
\includegraphics[scale=1]{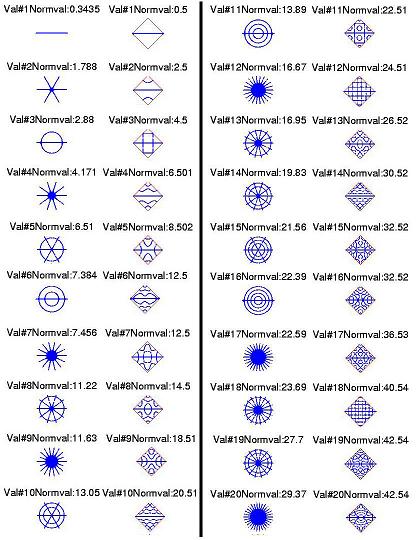}
\end{center}
\caption{$(2)$ Eigenfunction correspondences}
\label{figeleven5}
\end{figure}

    Observe from Figs. \ref{figone}-\ref{figfive} that the endpoints of the graphs of the eigenvalues give a bijection between the eigenfunctions on the circle and those on the square.  Since this bijection may encapsulate information about our homotopy or the eigenfunctions themselves, we examine it.  For simplicity, we view the correspondences of the nodal sets of the eigenfunctions.  Figs. \ref{figeleven1}-\ref{figeleven5} displays pictorially some of these correspondences and Table \ref{tableone} displays all of the data in terms of eigenfunction labels.  For instance, we find nontrivial information about the splitting of the two-dimensional eigenspaces on the endpoints (corresponding to $t=0$ and $t=1$).

How do we actually determine existence or non-existence of a two-dimensional eigenspace as we vary $t$ in an eigenfunction homotopy?  Since different symmetry families do not interact, it is easy to see when two eigenvalues from different symmetry families cross.  For instance, if we combine the (non-interacting) eigenvalues of Figs. \ref{figone}-\ref{figfive} in a single plot, we will definitely have crossings.  But for eigenfunctions in the same symmetry family, we cannot rigorously determine from the experimental data whether or not two eigenvalues cross, since we are only approximating the smooth functions $F$ and $G$ in discrete time.  Yet, as we vary the discrete time parameter, the graphs of the (lower eigenvalue) eigenfunctions change very slowly.  So, if two eigenvalues are far apart, and their corresponding eigenfunctions look similar as we vary $t$, then (trivially) these two eigenvalues do not attain the same value.  For example, in Fig. \ref{figone}, eigenvalues $1$ and $10$ on the circle (where $t=0$) are far apart, and their eigenfunctions look entirely different as we vary $t$.  Therefore, these eigenvalues do not cross.

Subtlety arises when two eigenvalues are very close together.  If a multiple eigenvalue exists, our numerical methods will not give us two identical eigenvalues.  This effect is a result of the instability of two-dimensional eigenspaces discussed in Section \ref{secintro}.  Since we are given approximate solutions to the eigenvalue problem, we will get two distinct eigenvalues which are very close together.  Measuring the difference will tell us the approximate (relative) error of our numerical method.  Given two nearly identical computed eigenvalues $\tilde{\lambda}_{j},\tilde{\lambda}_{j+1}$ at $t=t_{0}$, we are therefore interested in the quantity
\begin{equation}\label{two3}
E_{j}=\frac{\abs{\tilde{\lambda}_{j}-\tilde{\lambda}_{j+1}}}{\tilde{\lambda}_{j}}
\end{equation}

For eigenspaces that we know have multiplicity, we find that $E_{j}\leq 8\times10^{-5}$.  Recall from above that, with respect to our eigenvalue data, a collision (see Fig. \ref{figthirteen.5}(a)) and a multiplicity are effectively indistinguishable for $E_{j}$ very small.  Therefore, for collisions with $E_{j}$ around $8\times10^{-5}$, we are forced to conclude that a multiple eigenvalue occurs.  Conversely, if $E_{j}\gg 8\times10^{-5}$, we know that no multiple eigenvalue occurs in a neighborhood of the region we have analyzed.

In Fig. \ref{figone}, eigenvalues $12$ and $13$ on the circle (where $t=0$) become very close around $t=.4$ and $t=.5$.  We have observed the discrete homotopy of eigenvalue $12$ in Figs. \ref{figsix}, \ref{figseven} and \ref{figeight}.  In all three figures, we see what appears to be a discontinuity in the homotopy around $t=.5$.

This apparent discontinuity occurs since the numerical computations we have used cannot tell the difference between a two-dimensional eigenspace and an almost two dimensional eigenspace.  In other words, for two close eigenvalues, the solver combines two distinct eigenfunctions.  Our intuition used above can also apply here:  quite possibly, an eigenvalue collision is just a perturbation (geometric or otherwise) of a true two-dimensional eigenspace.  Once again, this is a manifestation of the superposition effect \cite{backer97}.  The data shows two linear combinations of the true eigenfunctions, instead of a separation into two non-interacting eigenspaces.  To assuage these concerns, we make higher precision calculations with a higher number of time values $t$ around the potential two-dimensional eigenspace.  We then use the relative error of Eq. (\ref{two3}) to determine whether or not a two-dimensional eigenspace occurs.  If a two-dimensional eigenspace is determined to occur, we then incorporate the eigenfunction data and the fact that the eigenvalues vary analytically (recall Theorem 3.9,VII.\S4.5 of \cite{kato66}).  If the higher and lower eigenfunctions swap appearances around the multiple eigenvalue, then the eigenvalues cross.  If no such swap occurs, the eigenvalues do not cross.  This completes our analysis of the spectral dynamics.

\begin{figure}[htbp!]
\begin{center}
\includegraphics[scale=1]{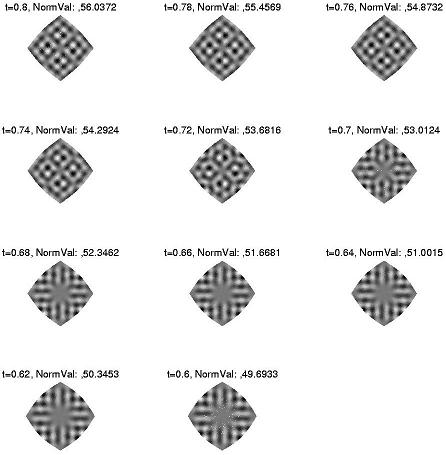}
\end{center}
\caption{$(1--)$ Family Tenth Eigenvalue Zoom}
\label{figtwelve}
\end{figure}

\begin{figure}[htbp!]
\begin{center}
\includegraphics[scale=1]{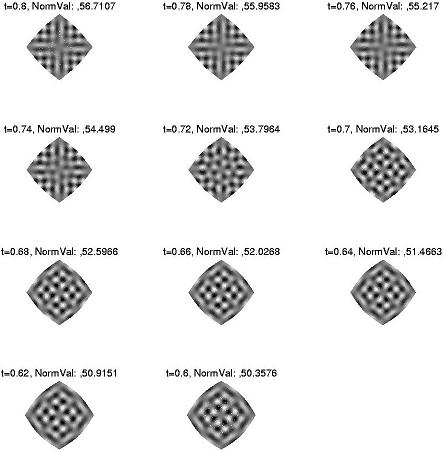}
\end{center}
\caption{$(1--)$ Family Eleventh Eigenvalue Zoom}
\label{figthirteen}
\end{figure}

\begin{figure}[htbp!]
\begin{center}
\includegraphics[scale=1]{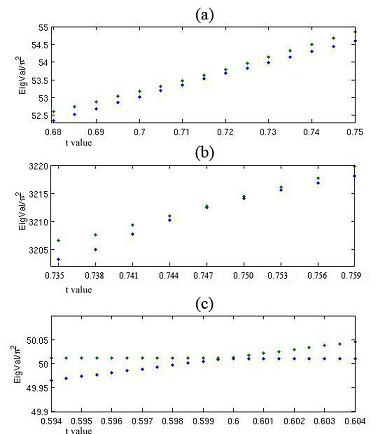}
\end{center}
\caption{Eigenvalue Dynamics Summary.
(a) A \textit{collision}.  Tenth and eleventh $(1--)$ eigenvalues of $H(t,\cdot,\cdot)$ as in Figs. \ref{figtwelve} and \ref{figthirteen}.
(b) A \textit{non-crossing} tangency, with (analytically determined) two-dimensional eigenspace.  These are two $(1++)$ eigenvalues of the homotopy $G_{0}(t,\cdot,\cdot)$ of Section \ref{seccarp} as in Figs. \ref{figtwentysix} and \ref{figtwentysix.1}.
(c) A \textit{crossing}, with (experimentally determined) two-dimensional eigenspace.  These are the sixth and seventh $(1+-)$ eigenvalues of the homotopy $G_{0}(t,\cdot,\cdot)$ of Section \ref{seccarp} as in Figs. \ref{figtwentyfour} and \ref{figtwentyfour.1}.}
\label{figthirteen.5}
\end{figure}

Let us demonstrate the above procedure by example.  Along the homotopy $H$ for the $(1--)$ family, we see from Fig. \ref{figfour} that eigenvalues $10$ and $11$ on the circle ($t=0$) get very close together at $t=.7$.  Our depiction of the homotopy of Eigenvalue $10$ on the circle in Fig. $\ref{fignine}$ again shows ambiguous behavior at $t=.7$.  To make our analysis more precise, we solve the eigenvalue problem for $H$ for $t\in\{.6,.62,.64,\ldots,.78,.8\}$.  In Fig. $\ref{figtwelve}$ we show the tenth eigenvalue for these $t$ values, and in Fig. $\ref{figthirteen}$ the eleventh eigenvalue.  In the time interval $[.6,.8]$, we observe that the tenth and eleventh eigenvalues become very close.  For instance, at $t=.72$, the normalized eigenvalues are separated by about $1/10$, giving a relative error $E_{10}\approx .00185\gg8\times10^{-5}$.

\begin{figure}[htbp!]
    \begin{center}
    \includegraphics[scale=1]{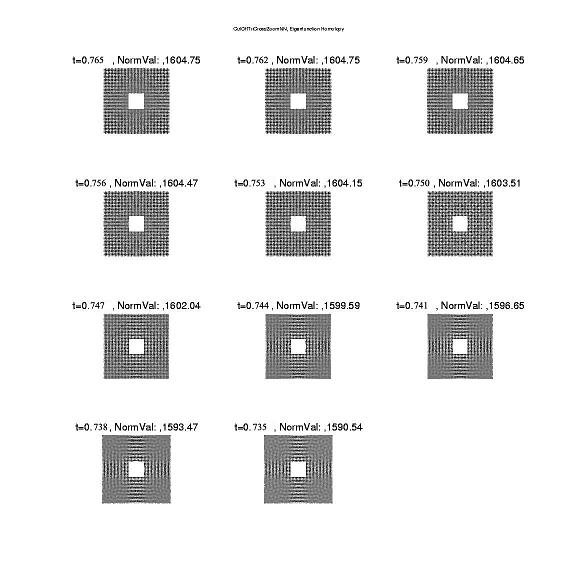}
    \end{center}
    \caption{$\lambda=1600$ at $t=.125$, EigFcn 1}
    \label{figtwentysix}
    \end{figure}

\begin{figure}[htbp!]
    \begin{center}
    \includegraphics[scale=1]{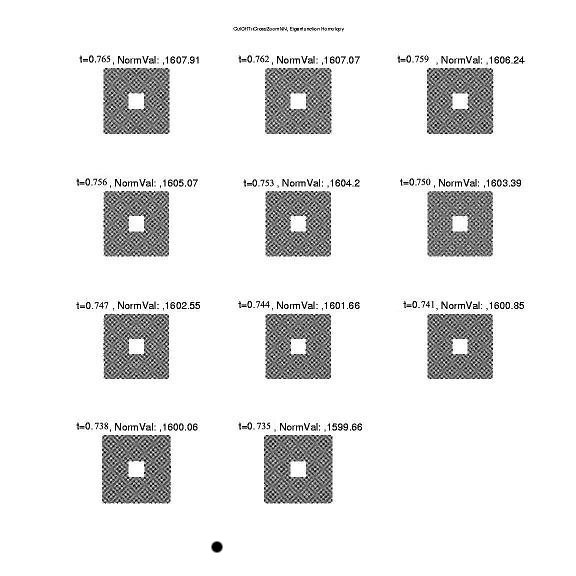}
    \end{center}
    \caption{$\lambda=1600$ at $t=.125$, EigFcn 2}
    \label{figtwentysix.1}
    \end{figure}

So, we conclude that no two-dimensional eigenspace has occurred within this time interval.  For an example with a two-dimensional eigenspace with no apparent crossing (i.e. where the ordering of the two eigenvalues does not change), see Figs \ref{figtwentysix} and \ref{figtwentysix.1} in Section \ref{seccarp}.

\begin{figure}[htbp!]
\begin{center}
\includegraphics[scale=1]{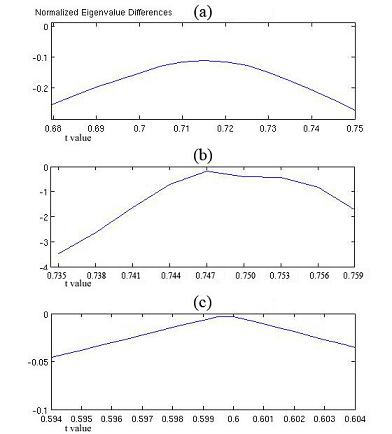}
\end{center}
\caption{Eigenvalue Dynamics Summary, Eigenvalue Differences $\lambda_{j+1}-\lambda_{j}$ from Fig. \ref{figthirteen.5}.
(a) A \textit{collision}.  Tenth and eleventh $(1--)$ eigenvalues of $H(t,\cdot,\cdot)$ as in Figs. \ref{figtwelve} and \ref{figthirteen}.
(b) A \textit{non-crossing} tangency, with (analytically determined) two-dimensional eigenspace.  These are two $(1++)$ eigenvalues of the homotopy $G_{0}(t,\cdot,\cdot)$ of Section \ref{seccarp} as in Figs. \ref{figtwentysix} and \ref{figtwentysix.1}.
(c) A \textit{crossing}, with (experimentally determined) two-dimensional eigenspace.  These are the sixth and seventh $(1+-)$ eigenvalues of the homotopy $G_{0}(t,\cdot,\cdot)$ of Section \ref{seccarp} as in Figs. \ref{figtwentyfour} and \ref{figtwentyfour.1}.}
\label{figthirteen.6}
\end{figure}

\begin{figure}[htbp!]
    \begin{center}
    \includegraphics[scale=1]{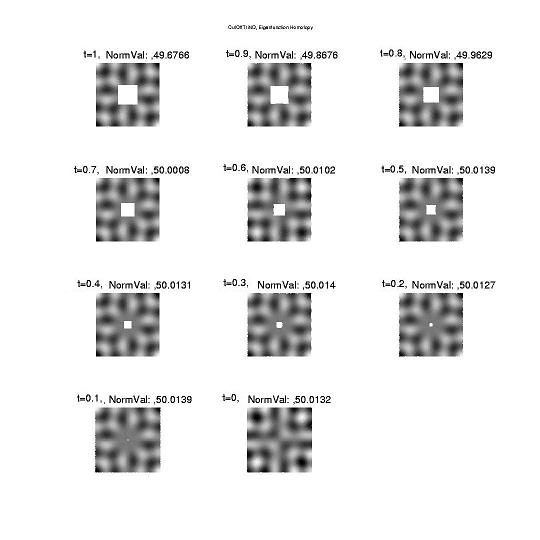}
    \end{center}
    \caption{$(1+-)$ EigVal 6}
    \label{figtwentyfour}
    \end{figure}

\begin{figure}[htbp!]
    \begin{center}
    \includegraphics[scale=1]{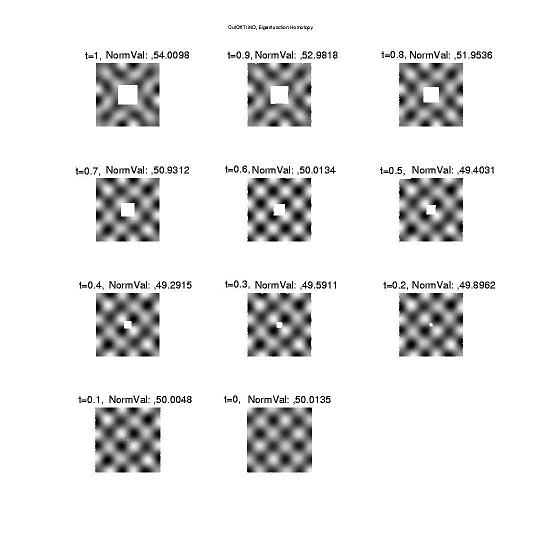}
    \end{center}
    \caption{$(1+-)$ EigVal 7}
    \label{figtwentyfour.1}
    \end{figure}

To summarize the different eigenvalue dynamics, we display Figs. \ref{figthirteen.5} and \ref{figthirteen.6}.  Fig. \ref{figthirteen.5} shows three different eigenvalue interactions (all from eigenvalue computations), and Fig. \ref{figthirteen.6} shows the gaps between each of the respective pairs of eigenvalues.  If we examine the axes of these figures carefully, we see that part (a) has a pair of eigenvalues which do not come close enough together for a two-dimensional eigenspace to occur.  This reaffirms our recent observations of Figs. \ref{figtwelve} and \ref{figthirteen}.  From these same figures we see that the two eigenfunctions swap identities, so with our intuition from Eq. (\ref{one2.5}) in mind, we call this interaction a \textit{collision}.  In part (b), we know that a two-dimensional eigenspace occurs at $t=.75$, but we see from Figs. \ref{figtwentysix} and \ref{figtwentysix.1} that the eigenfunctions do not swap identities.  Therefore, this is a \textit{non-crossing}.  In part (c), we find evidence for a two-dimensional eigenspace, and we observe the swapping of eigenfunction identities in Figs. \ref{figtwentyfour} and \ref{figtwentyfour.1}.  Therefore, this is a \textit{crossing}.

We note that the same effect of nearby eigenspaces also occurs on the square.  Recall that on the square, we have several two dimensional eigenspaces.  So, the choice of two orthogonal eigenfunctions by our numerical method from such an eigenspace amounts to two arbitrary linear combinations of the true eigenfunctions.  However, as in Fig. \ref{figsix}, we see that, as we approach the square (around $t=.9$), we approach a particular eigenfunction in a two-dimensional eigenspace on the square.  Since the eigenfunctions vary analytically, the homotopy of eigenfunctions therefore selects for us a member of the two dimensional eigenspace on the square.  It is then interesting to note that, in our observed data, the eigenfunction thus selected is always a ``natural'' eigenfunction on the square, i.e. an eigenfunction of the form $\cos\pi ax\cos\pi by\pm\cos\pi bx\cos\pi ay$, rather than a linear combination of two such eigenfunctions.  However, in the homotopy of Section \ref{seccarp}, we often see that a two-dimensional eigenspace on the square splits into two one-dimensional eigenspaces in the non-``natural'' way.  More specifically, the numerical evidence shows a few two-dimensional eigenspaces definitely split in a non-canonical way (see Figs. \ref{figtwentyfive} and \ref{figtwentyfive.1} below).  For other examples, the evidence is less striking, but non-canonical splittings still seem to occur in mostly all examples (bold entries in Table \ref{tabletwo}).

\section{Sierpinski Carpet Homotopy}
\label{seccarp}

Let $A_{0}=\left\{(x,y)\in\R^{2}\colon\max\{{\abs{x},\abs{y}}\}\leq\frac{1}{2}\right\}$ be the unit square centered at the origin.  Let $\{f_{i}\}_{i=1}^{8}$ denote the usual eight affine linear maps used in the construction of $SC$, and let
\begin{equation}\label{three1}
A_{j}=\bigcup_{i=1}^{8}f_{i}A_{j-1},
\end{equation}
for $j>0$, $j\in\Z$.  Recall that the Sierpinski carpet $SC$ is defined (as in Section \ref{secintro}) by
\begin{equation}\label{three2}
SC=\bigcap_{i=0}^{\infty}A_{i}.
\end{equation}
We would like to deform $A_{1}$ continuously into $A_{0}$, and then iterate this procedure for higher $A_{j}$.  (The choice to make a map $A_{1}\to A_{0}$ instead of $A_{0}\to A_{1}$ is for convenience.  Also, as stated in Section \ref{secintro}, note that the removal of a single point does not change the Neumann Laplacian.)  To construct these deformations, we follow the procedure of Section \ref{secctos}.  Let
\begin{equation}\label{three2.5}
C=\left\{(x,y)\in\R^{2}\colon x>0,\frac{\pi}{4}>\tan^{-1}\left(y/x\right)>-\frac{\pi}{4}\right\}\cup\{(0,0)\}
\end{equation}
so that $C$ is the right ``quarter'' of the plane.  Then let $L=A_{1}\cap C$ be the right ``quarter'' of the picture frame $A_{1}$.  Now, define $G(t,x,y)\colon[0,1]\times L\to A_{0}\cap C$ by
\begin{equation}\label{three3}
G(t,x,y)=\left(\frac{3}{2}tx-\frac{3}{2}t+1\right)\cdot(x,y).
\end{equation}

\begin{figure}[htbp!]
    \begin{center}
    \includegraphics[scale=1]{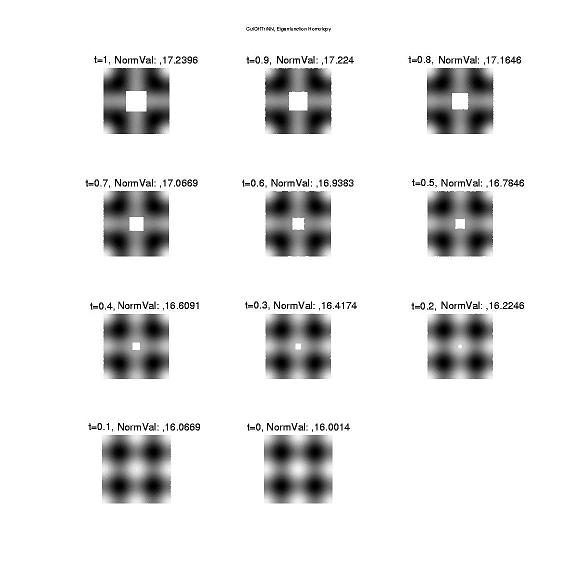}
    \end{center}
    \caption{$G_{0}, (1++)$ EigVal $4$ Homotopy}
    \label{fignineteen}
    \end{figure}

\begin{figure}[htbp!]
    \begin{center}
    \includegraphics[scale=1]{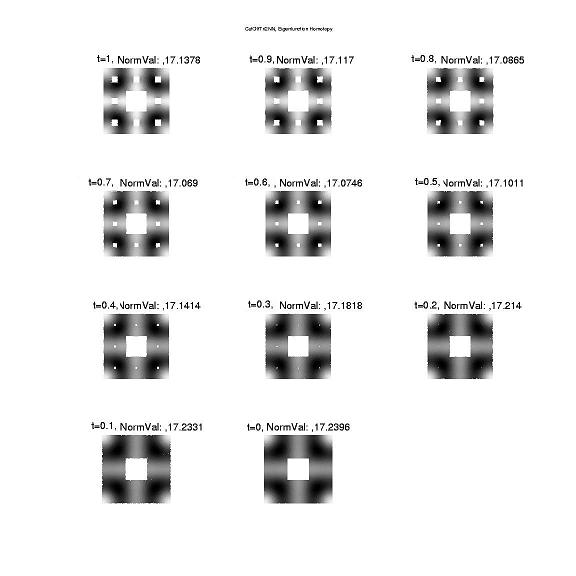}
    \end{center}
    \caption{$G_{1}, (1++)$ EigVal $4$ Homotopy}
    \label{figtwenty}
    \end{figure}

\begin{figure}[htbp!]
    \begin{center}
    \includegraphics[scale=1]{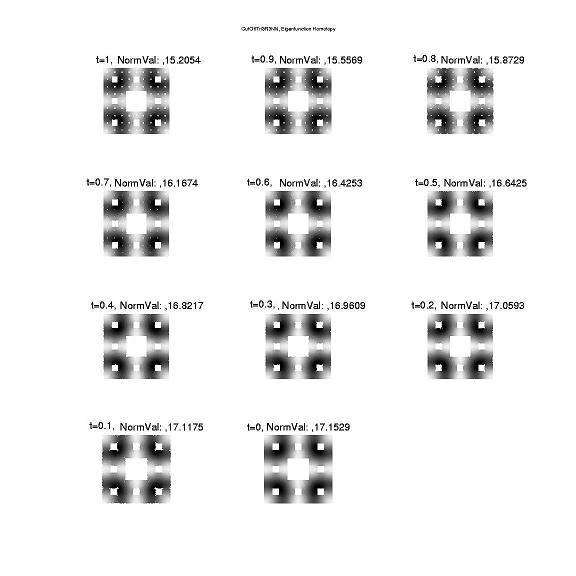}
    \end{center}
    \caption{$G_{2}, (1++)$ EigVal $4$ Homotopy}
    \label{figtwentyone}
    \end{figure}

\begin{figure}[htbp!]
    \begin{center}
    \includegraphics[scale=1]{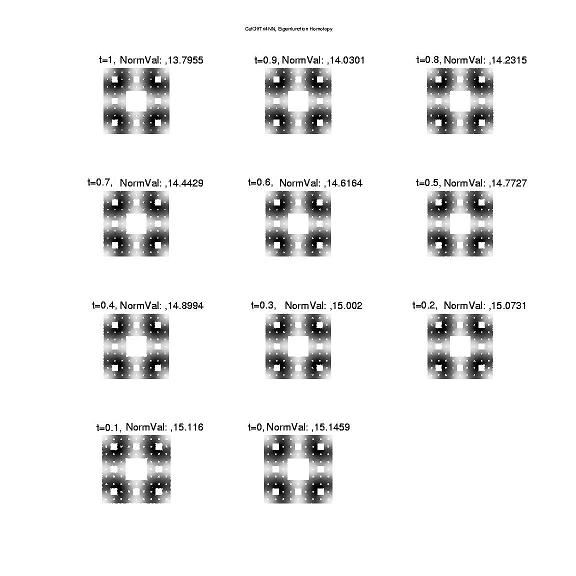}
    \end{center}
    \caption{$G_{3}, (1++)$ EigVal $4$ Homotopy}
    \label{figtwentytwo}
    \end{figure}

Observe that $G$ is a straight line homotopy from the identity map to the map $(x,y)\mapsto((3/2)x-(1/2))\cdot(x,y)$.  As in Section \ref{secctos}, if we rotate the domain and range of the map three times by increments of $\pi/2$ radians, we obtain a continuous map $G_{0}\colon [0,1]\times A_{1}\to A_{0}$.  This map is depicted in Fig. \ref{fignineteen}.  To extend $G_{0}$ to the remaining $A_{j}$, $j>0$, define $G_{j}(t,x,y)\colon [0,1]\times A_{j+1}\to A_{j}$ inductively by
\begin{equation}\label{three4}
G_{j}(t,x,y)=\bigcup_{i=1}^{8}f_{i}\circ G_{j-1}(t,f_{i}^{-1}(x,y)).
\end{equation}
In words, we apply $G_{j-1}$ to each of the eight images $\{f_{i}A_{j-1}\}_{i=1}^{8}$ used to define $A_{j}$.  So, for a fixed $j$, we use $f_{i}\colon A_{j-1}\to A_{j}$, so that $f_{i}^{-1}\colon f_{i}A_{j-1}\to A_{j-1}$ is the (well-defined) map in Eq. (\ref{three4}).  Finally, since $G_{0}=id$ on $\partial A_{0}\subset A_{1}$ for all $t$, the map $G_{j}$ is well defined (a priori, each $G_{j}(t,\cdot,\cdot)$ has two definitions for certain $(x,y)\in B_{j}$).  The maps $G_{1}$ through $G_{3}$ are depicted in Figs. \ref{figtwenty} through \ref{figtwentytwo} respectively.

As in Section \ref{secctos} we see that, for each $j$, the eigenvalues and eigenvectors of $G_{j}$ can be represented as functions holomorphic in $t$.  However, since the domains of $G_{j}$ and $G_{j+1}$ only meet continuously, the eigenvalues only vary continuously around $\im \left(G_{j}(1,\cdot,\cdot)\right)=\im \left(G_{j+1}(0,\cdot,\cdot)\right)$.  We further note that $G_{j}(1,x,y)$ maps $B_{j}$ to $A_{j+1}$ as expected. Finally, to stay consistent with Section \ref{secctos}, we refer to $G_{j}(t,x,y)$ below as $G_{j}(1-t,x,y)$ as defined in the previous few paragraphs (we have already used this convention in the caption of Figs. \ref{figthirteen.5} and \ref{figthirteen.6}).

As above, our domains $\{B_{j}\}_{j=0}^{\infty}$ are invariant under the $D_{4}$ symmetry group.  In our computations, we restrict each $B_{j}$ to the sector $0\leq\theta\leq\pi/4$ in polar coordinates, for the $1$-dimensional representations.  We impose our four combinations of boundary conditions on the intersection of our domain with the rays $\theta=0$ and $\theta=\pi/4$.  For the $2$-dimensional representation, our fundamental domain is $B_{j}\cap L$.  We impose Neumann boundary conditions on the ray $\theta=\pi/4$ and Dirichlet conditions on the ray $\theta=-\pi/4$.

For all fundamental domains, we place Neumann boundary conditions on all remaining edges.  Notice that once again, our fundamental domains are simply connected for $j=0$.  Therefore, as in Theorem 6.4 of \cite{teytel99}, we should expect no eigenvalue crossings along this homotopy (for each individual symmetry family).  Since eigenvalue crossings do in fact occur, this suggests that there is something exceptional (i.e. non-generic) about our domains.  The appearance of these crossings could be related to the quantum-classical correspondence principle.  More specifically, since the billiard flow is not ergodic in phase space, this may suggest a dense number of crossings.

\begin{figure}[htbp!]
    \begin{center}
    \includegraphics[scale=1]{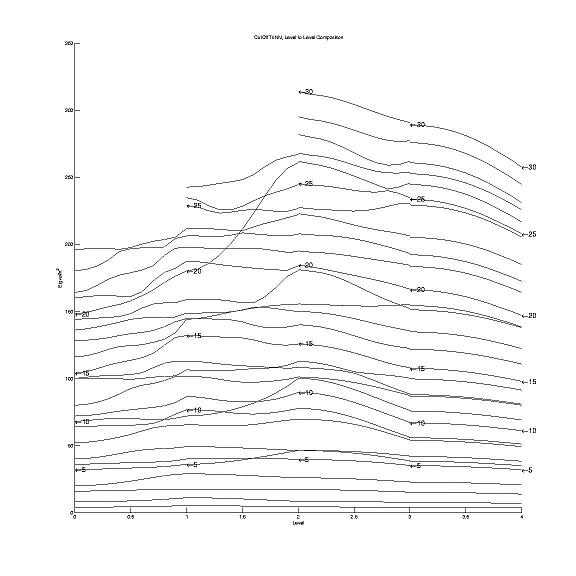}
    \end{center}
    \caption{$G_{0}$ through $G_{3}, (1++)$ Eigenvalue Summary Plot}
    \label{figfourteen}
    \end{figure}

\begin{figure}[htbp!]
    \begin{center}
    \includegraphics[scale=1]{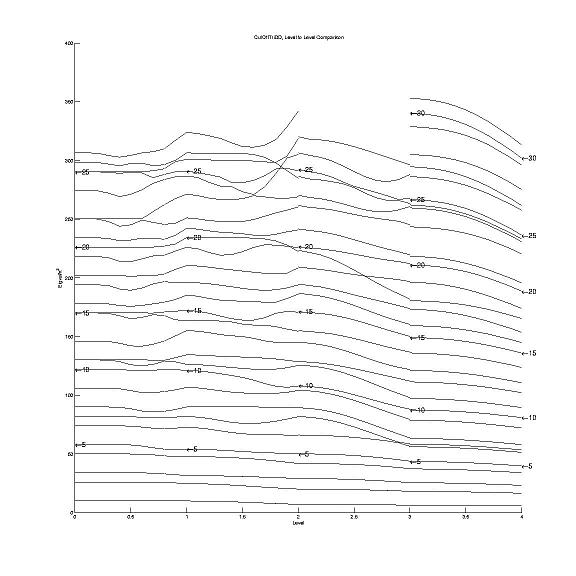}
    \end{center}
    \caption{$G_{0}$ through $G_{3}, (1--)$ Eigenvalue Summary Plot}
    \label{figfifteen}
    \end{figure}

\begin{figure}[htbp!]
    \begin{center}
    \includegraphics[scale=1]{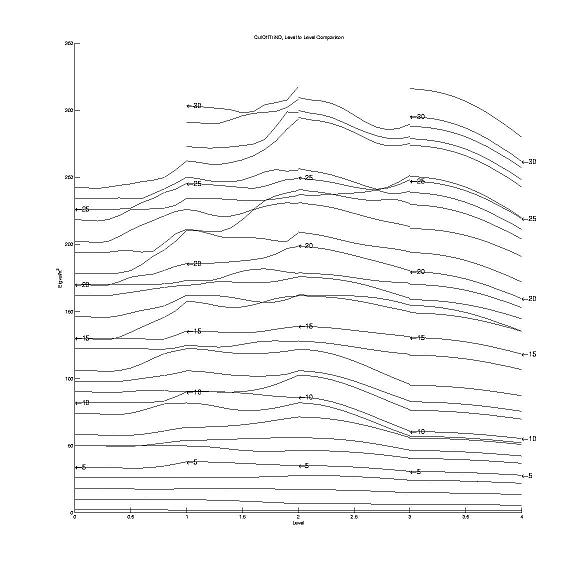}
    \end{center}
    \caption{$G_{0}$ through $G_{3}, (1-+)$ Eigenvalue Summary Plot}
    \label{figsixteen}
    \end{figure}

\begin{figure}[htbp!]
    \begin{center}
    \includegraphics[scale=1]{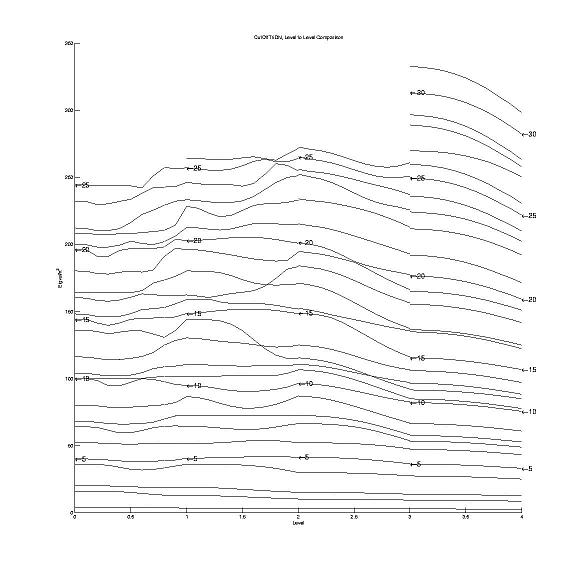}
    \end{center}
    \caption{$G_{0}$ through $G_{3}, (1+-)$ Eigenvalue Summary Plot}
    \label{figseventeen}
    \end{figure}

\begin{figure}[htbp!]
    \begin{center}
    \includegraphics[scale=1]{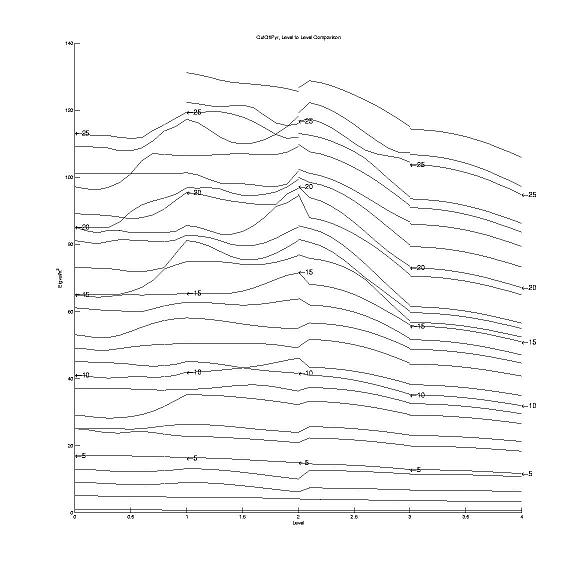}
    \end{center}
    \caption{$G_{0}$ through $G_{3}, (2)$ Eigenvalue Summary Plot}
    \label{figeighteen}
    \end{figure}

    In Figs. \ref{figfourteen} to \ref{figeighteen} we display the observed crossing of eigenvalues for the five symmetry families ($1++,1+-,1-+,1--,2$) with eigenvalues divided by $\pi^{2}$.  We collect the data of $G_{j}(t,\cdot,\cdot)$ for $j=0,1,2,3$ and for $t\in\{k/10\}_{k=0}^{10}$.  Observe that, in the bottom of the spectrum , the incidence of crossings decreases as $j$ increases.  More specifically, the incidence of non-crossings grows exponentially.  Also observe that, indeed, the eigenvalues of $G_{j}$ and $G_{j+1}$ meet in a continuous but non-differentiable manner.

    In Figs. \ref{fignineteen}-\ref{figtwentytwo}, we show the homotopy of eigenfunctions of Eigenvalue $4$ in the $(1++)$ family for $G_{0}$ through $G_{3}$.  This eigenfunction exhibits the typical, expected behavior.  The function changes a great deal for $G_{j}$ with $j$ small, and then as $j$ grows large the eigenfunction maintains the same appearance.  However, the eigenvalue continues changing, as expected.

\begin{figure}[htbp!]
    \begin{center}
    \includegraphics[scale=1]{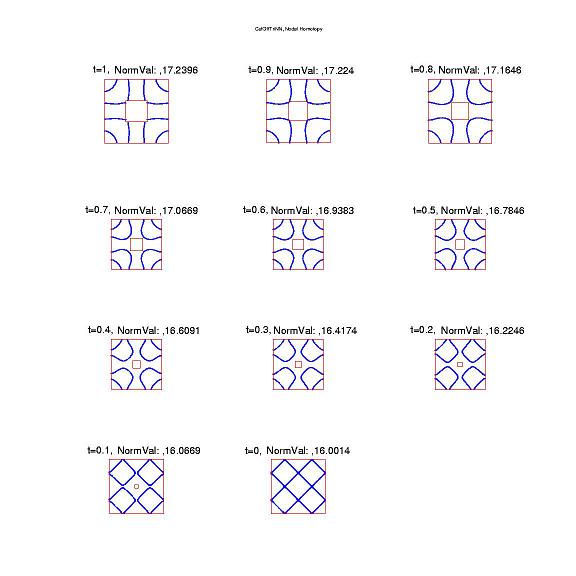}
    \end{center}
    \caption{$G_{0}, (1++)$ EigVal $4$ Nodal Homotopy}
    \label{figtwentythree}
    \end{figure}

\begin{figure}[htbp!]
    \begin{center}
    \includegraphics[scale=1]{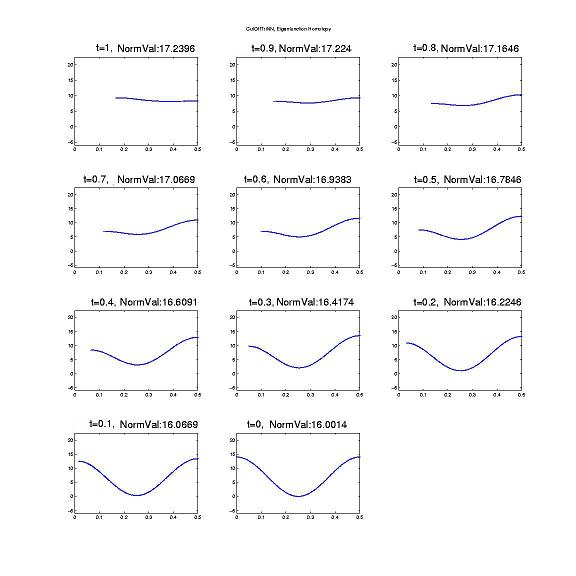}
    \end{center}
    \caption{$G_{0}, (1++)$ EigVal $4$ Line Homotopy}
    \label{figtwentythree.1}
    \end{figure}

    In Fig. \ref{figtwentythree} we show the nodal homotopy corresponding to Fig. \ref{fignineteen}.  With a low eigenvalue, this eigenfunction exhibits minimal irregularity during the homotopy.  In Fig. \ref{figtwentythree.1} we show the homotopy of \ref{fignineteen} restricted to the line $\{(x,y)\in\R^{2}\colon 0\leq x\leq1/2\}$.

    It is interesting to note how the two dimensional eigenspaces split for $G_{0}(t,\cdot,\cdot)$ as we move $t$ from $0$ to some small $\epsilon>0$.  In Section \ref{secctos}, we saw a ``canonical splitting'' for the given homotopy.  However, the eigenfunctions for $G_{0}(\epsilon,\cdot,\cdot)$ for two dimensional eigenspaces do not resemble the natural eigenfunctions of $G_{0}(0,\cdot,\cdot)$.  For example, in Figs. \ref{figtwentyfour} and \ref{figtwentyfour.1}, we see the two-dimensional eigenspace on the square split with an apparent discontinuity.  In Figs. \ref{figtwentyfive} and \ref{figtwentyfive.1}, we see a (non-canonical) splitting of the two-dimensional eigenspace on the square.  If this were a canonical splitting, we would have the $(10,0)$ and $(8,6)$ eigenfunctions.  However, there is no direct resemblence of our computed eigenfunctions and the canonical ones.  In fact, around $t=.1$, these eigenfunctions cross again, since the relative error between the eigenvalues is about $6.4\times10^{-5}<8\times10^{-5}$ (recall the discussion around Eq. (\ref{two3}).  Finally, in Figs. \ref{figtwentysix} and \ref{figtwentysix.1} we show the homotopy around an (analytically determined) two-dimensional eigenspace at $t=.75$.  It is not hard to find infinite sets of such eigenspaces for any rational $t\in(0,1)$.

\section{Conclusions}
\label{secconc}

We now summarize our results.  In both of our homotopies, we preserve the $D_{4}$ symmetry of the domains.  This allows us to split the eigenfunctions on these domains into five non-interacting families.  That is, the behavior of eigenfunctions in one symmetry family along the homotopy does not affect the other families.  The following is then clear: eigenvalues of one symmetry family definitely cross with eigenvalues of other symmetry families.  However, the numerical evidence leads us to speculate on the individual symmetry families as follows.  When we deform the circle to the square, we observe no intermediate two-dimensional eigenspaces.  However, when we deform the square to the picture frame (a square with a smaller square removed from the center) we do find intermediate two-dimensional eigenspaces.  In fact, analytically one can easily find an arbitrarily large number nontrivial multiplicity eigenspaces for a fixed rational time.

When we observe our domains, we see that the deformation of the circle to the square is a convex diamond-shaped region for all $t\in[0,1]$.  In fact, for $t\in(0,1)$, the region is a diamond with four edges, each of which is a half-ellipse.  However, as we map the square $A_{0}$ to the picture frame $A_{1}$ via $H$, we always have a rational polygon (a polygon with angles which are rational multiples of $\pi$).  Mapping $A_{i}$ to $A_{i+1}$ via $G_{i}$ always yields a rational domain with a finite number of flat edges (which is now multiply connected for $i\geq1$).  Recall: (a) after reducing multiplicities due to symmetry, we observe eigenvalues of nontrivial multiplicity for $G_{0}$ but not $H$.  Also, (b) for $G_{i}$ we observe non-canonical splitting of two-dimensional eigenspaces of the square, and for $H$ we only observe canonical splitting.

With respect to observation (a), it is tempting to make a connection between the appearance of eigenvalue multiplicity and ergodicity of the billiard flow, but such a connection is unclear.  Such considerations seem related to the quantum-classical correspondence, and some conjectures discussed in \cite{sarnak92}.  In particular, the author mentions that, for the bottom part of the spectrum, the eigenvalues for integrable billiards may follow a Poisson distribution, while the eigenvalues for a chaotic system may follow a GOE (Gaussian Orthogonal Ensemble) distribution.  With respect to observation (b), the authors have no clear explanation.

We have just discussed possible connections between ergodicity of billiard flow and eigenvalue properties.  Such connections are collected under the so-called quantum-classical correspondence.  Perhaps extending these connections further, we discuss the localization and near-localization (or scarring) of eigenfunctions.  For the standard Laplacian on the Sierpinski Gasket ($SG$), many eigenfunctions are localized \cite{strichartz06}.  Indeed, it appears that for many symmetric PCF fractals, the symmetry of the fractal and the capability to break it into disjoint pieces by removing finitely many points contributes to the localization of eigenfunctions. In contrast, we have $SC$, which is a symmetric non-PCF fractal (intuitively, this means it cannot be broken up into smaller pieces by removing finitely many points).

We now ask: do the ``stronger'' connectivity properties of $SC$ forbid localized eigenfunctions?  We believe the numerical evidence from this paper suggests, yes.  Indeed, no localized eigenfunctions appear in the computations.  If the outer approximation method of \cite{berry08} actually works, then this observation finds strong support.  For possible counter-evidence, let us consider the results of \cite{hassell09}, which imply that any rectangular portion of a domain almost always exhibits eigenfunction scarring (that is, having mass concentrated in the domain in a non-equidistributed way).  This would imply that the number of scarred modes for our domains should increase exponentially as $j\to\infty$ as in Eq. (\ref{one10}).  However, it is unclear whether a scarred eigenfunction for a fixed $j$ would remain scarred as $j$ increases.  If the outer approximation method of \cite{berry08} holds, then the eigenfunction scarring may have repercussions for the eigenfunctions on $SC$, even though strict localization may not occur.


Although it may seem unfortunate that we are not able to arrive at any coherent description, even just conjectural, relating the behavior of eigenfunctions joined by a homotopy, there is also a more optimistic conclusion: the connection is rather complicated, and requires further investigation.  We hope that the data we present will inspire further investigation.  (Data available in full at \verb!www.math.cornell.edu/~smh82!)


\pgbreak

\begin{table}[htbp!]
\resizebox{13cm}{!}{
\begin{tabular}{r||rr||rr||rr||rr||rr}
\multicolumn{10}{p{9cm}}{Eigenvalue Correspondences}\\
\hline
\multicolumn{1}{p{2cm}}{Eigenvalue Number on Circle}
& \multicolumn{2}{p{2.5cm}}{$(1++)$ Family}
& \multicolumn{2}{p{2.5cm}}{$(1--)$ Family}
& \multicolumn{2}{p{2.5cm}}{$(1-+)$ Family}
& \multicolumn{2}{p{2.5cm}}{$(1+-)$ Family}
& \multicolumn{2}{p{2.5cm}}{$(2)$ Family}\\
\hline
& Circle & Square & Circle & Square & Circle & Square & Circle & Square & Circle & Square\\
\hline
  1 & (      1,      0) & (      2,      0)          & (      0,      4) & (      3,      1)          & (      0,      2) & (      1,      1)          & (      0,      2) & (      2,      0)          & (      0,      1) & (      1,      0)\\
  2 & (      0,      4) & (      2,      2)          & (      1,      4) & (      5,      1)          & (      1,      2) & (      3,      1)          & (      1,      2) & (      4,      0)          & (      0,      3) & (      2,      1)\\
  3 & (      2,      0) & (      4,      0)          & (      0,      8) & (      5,      3)          & (      0,      6) & (      3,      3)          & (      0,      6) & (      4,      2)          & (      1,      1) & (      3,      0)\\
  4 & (      1,      4) & (      4,      2)          & (      2,      4) & (      7,      1)          & (      2,      2) & (      5,      1)          & (      2,      2) & (      6,      0)          & (      0,      5) & (      3,      2)\\
  5 & (      0,      8) & (      4,      4)          & (      0,     12) & (      7,      3)          & (      1,      6) & (      5,      3)          & (      1,      6) & (      6,      2)          & (      1,      3) & (      4,      1)\\
  6 & (      3,      0) & (      6,      0)          & (      1,      8) & (      7,      5)          & (      0,     10) & \textbf{(      5,      5)} & (      0,     10) & (      6,      4)          & (      2,      1) & \textbf{(      4,      3)}\\
  7 & (      2,      4) & (      6,      2)          & (      3,      4) & (      9,      1)          & (      3,      2) & \textbf{(      7,      1)} & (      3,      2) & (      8,      0)          & (      0,      7) & \textbf{(      5,      0)}\\
  8 & (      4,      0) & (      6,      4)          & (      2,      8) & (      9,      3)          & (      2,      6) & (      7,      3)          & (      2,      6) & (      8,      2)          & (      1,      5) & (      5,      2)\\
  9 & (      0,     12) & (      8,      0)          & (      0,     16) & (      9,      5)          & (      0,     14) & (      7,      5)          & (      0,     14) & (      8,      4)          & (      0,      9) & (      6,      1)\\
 10 & (      1,      8) & (      8,      2)          & (      1,     12) & (     11,      1)          & (      4,      2) & (      9,      1)          & (      4,      2) & \textbf{(      8,      6)} & (      2,      3) & (      5,      4)\\
 11 & (      3,      4) & (      6,      6)          & (      4,      4) & \textbf{(      9,      7)} & (      1,     10) & (      9,      3)          & (      1,     10) & \textbf{(     10,      0)} & (      3,      1) & (      6,      3)\\
 12 & (      5,      0) & (      8,      4)          & (      3,      8) & \textbf{(     11,      3)} & (      3,      6) & (      7,      7)          & (      3,      6) & (     10,      2)          & (      0,     11) & (      7,      0)\\
 13 & (      2,      8) & \textbf{(     10,      0)} & (      0,     20) & (     11,      5)          & (      5,      2) & (      9,      5)          & (      5,      2) & (     10,      4)          & (      1,      7) & (      7,      2)\\
 14 & (      0,     16) & \textbf{(      8,      6)} & (      5,      4) & \textbf{(     13,      1)} & (      0,     18) & (     11,      1)          & (      0,     18) & (     10,      6)          & (      2,      5) & (      6,      5)\\
 15 & (      1,     12) & (     10,      2)          & (      2,     12) & \textbf{(     11,      7)} & (      2,     10) & \textbf{(     11,      3)} & (      2,     10) & (     12,      0)          & (      3,      3) & \textbf{(      7,      4)}\\
 16 & (      4,      4) & (     10,      4)          & (      1,     16) & (     13,      3)          & (      1,     14) & \textbf{(      9,      7)} & (      1,     14) & (     12,      2)          & (      4,      1) & \textbf{(      8,      1)}\\
 17 & (      6,      0) & (      8,      8)          & (      4,      8) & (     13,      5)          & (      4,      6) & (     11,      5)          & (      4,      6) & (     12,      4)          & (      0,     13) & (      8,      3)\\
 18 & (      3,      8) & (     10,      6)          & (      6,      4) & (     11,      9)          & (      6,      2) & (      9,      9)          & (      6,      2) & (     10,      8)          & (      1,      9) & (      9,      0)\\
 19 & (      0,     20) & (     12,      0)          & (      3,     12) & (     13,      7)          & (      3,     10) & (     13,      1)          & (      3,     10) & (     12,      6)          & (      2,      7) & \textbf{(      9,      2)}\\
 20 & (      5,      4) & (     12,      2)          & (      0,     24) & (     15,      1)          & (      0,     22) & (     11,      7)          & (      0,     22) & (     14,      0)          & (      0,     15) & \textbf{(      7,      6)}\\
 21 & (      2,     12) & (     12,      4)          &                   &                            & (      2,     14) & (     13,      3)          & (      2,     14) & (     14,      2)          & (      3,      5) & (      8,      5)\\
 22 & (      7,      0) & (     10,      8)          &                   &                            & (      5,      6) & (     13,      5)          & (      5,      6) & (     12,      8)          & (      1,     11) & (      9,      4)\\
 23 & (      1,     16) & (     12,      6)          &                   &                            &                   &                            &                   &                            &                   &                  \\
 24 & (      4,      8) & (     14,      0)          &                   &                            &                   &                            &                   &                            &                   &                  \\
 25 & (      6,      4) &                            &                   &                            &                   &                            &                   &                            &                   &                  \\

\end{tabular}
}
\caption{
Table of Eigenfunction Correspondences:
For the Circle, $(a,b)\leftrightarrow J_{b}(\lambda_{n,a}r)\cos(2\pi b\theta)$.
For the Square, $(a,b)\leftrightarrow \cos(\pi ax)\cos(\pi by)\pm\cos(\pi ay)\cos(\pi bx)$,
($+$ in columns $1,3$, $-$ in columns $2,4$).
Two dimensional eigenspaces in bold.
}
\label{tableone}
\end{table}

\pgbreak

\begin{table}[htbp!]
\resizebox{16cm}{!}{
\begin{tabular}{r|rrrrr||r|rrrrr||r|rrrrr||r|rrrrr||r|rrrrr}
\multicolumn{30}{p{9cm}}{Eigenvalue Correspondences}\\
\hline
\multicolumn{6}{p{2.5cm}}{$(1++)$ Family}
& \multicolumn{6}{p{2.5cm}}{$(1--)$ Family}
& \multicolumn{6}{p{2.5cm}}{$(1+-)$ Family}
& \multicolumn{6}{p{2.5cm}}{$(1-+)$ Family}
& \multicolumn{6}{p{2.5cm}}{$(2)$ Family}\\
\hline
\multicolumn{1}{p{1.4cm}}{Eigfcn on Square} & $A_{0}$ & $A_{1}$ & $A_{2}$ & $A_{3}$ & $A_{4}$ &
\multicolumn{1}{p{1.4cm}}{Eigfcn on Square} & $A_{0}$ & $A_{1}$ & $A_{2}$ & $A_{3}$ & $A_{4}$ &
\multicolumn{1}{p{1.4cm}}{Eigfcn on Square} & $A_{0}$ & $A_{1}$ & $A_{2}$ & $A_{3}$ & $A_{4}$ &
\multicolumn{1}{p{1.4cm}}{Eigfcn on Square} & $A_{0}$ & $A_{1}$ & $A_{2}$ & $A_{3}$ & $A_{4}$ &
\multicolumn{1}{p{1.4cm}}{Eigfcn on Square} & $A_{0}$ & $A_{1}$ & $A_{2}$ & $A_{3}$ & $A_{4}$ \\
\hline
(      2,      0)          &       1 &       1 &       1 &       1 &       1 & (      3,      1)                   &       1 &       1 &       1 &       1 &       1 & (      2,      0)          &       1 &       1 &       1 &       1 &       1 & (      1,      1)          &       1 &       1 &       1 &       1 &       1 & (      1,      0)          &       1 &       1 &       1 &       1 &       1\\
(      2,      2)          &       2 &       2 &       2 &       2 &       2 & (      5,      1)                   &       2 &       2 &       2 &       2 &       2 & (      4,      0)          &       2 &       2 &       2 &       2 &       2 & (      3,      1)          &       2 &       2 &       2 &       2 &       2 & (      2,      1)          &       2 &       2 &       2 &       2 &       2\\
(      4,      0)          &       3 &       3 &       3 &       3 &       3 & (      5,      3)                   &       3 &       3 &       3 &       3 &       3 & (      4,      2)          &       3 &       3 &       3 &       3 &       3 & (      3,      3)          &       3 &       3 &       3 &       3 &       3 & (      3,      0)          &       3 &       3 &       3 &       3 &       3\\
(      4,      2)          &       4 &       4 &       4 &       4 &       4 & (      7,      1)                   &       4 &       4 &       4 &       4 &       4 & (      6,      0)          &       4 &       4 &       4 &       4 &       4 & (      5,      1)          &       4 &       4 &       4 &       4 &       4 & (      3,      2)          &       4 &       4 &       4 &       4 &       4\\
(      4,      4)          &       5 &       5 &       6 &       6 &       6 & (      7,      3)                   &       5 &       5 &       5 &       5 &       5 & (      6,      2)          &       5 &       5 &       5 &       5 &       5 & (      5,      3)          &       5 &       5 &       5 &       5 &       5 & (      4,      1)          &       5 &       5 &       5 &       5 &       5\\
(      6,      0)          &       6 &       6 &       5 &       5 &       5 & (      7,      5)                   &       6 &       6 &       6 &       6 &       6 & (      6,      4)          &       6 &       6 &       6 &       6 &       6 & \textbf{(      7,      1)} &       6 &       6 &       6 &       6 &       6 & \textbf{(      4,      3)} &       6 &       7 &       7 &       7 &       7\\
(      6,      2)          &       7 &       7 &       7 &       7 &       7 & (      9,      1)                   &       7 &       7 &       7 &       7 &       7 & (      8,      0)          &       7 &       7 &       7 &       7 &       7 & \textbf{(      5,      5)} &       7 &       7 &       7 &       7 &       7 & \textbf{(      5,      0)} &       7 &       6 &       6 &       6 &       6\\
(      6,      4)          &       8 &       8 &       8 &       8 &       8 & (      9,      3)                   &       8 &       8 &       8 &       8 &       8 & (      8,      2)          &       8 &       8 &       8 &       8 &       8 & (      7,      3)          &       8 &       8 &       8 &       8 &       8 & (      5,      2)          &       8 &       8 &       8 &       8 &       8\\
(      8,      0)          &       9 &       9 &      11 &      11 &      11 & (      9,      5)                   &       9 &       9 &       9 &       9 &       9 & (      8,      4)          &       9 &       9 &       9 &       9 &       9 & (      7,      5)          &       9 &       9 &       9 &       9 &       9 & (      6,      1)          &       9 &       9 &       9 &       9 &       9\\
(      8,      2)          &      10 &      10 &       9 &       9 &       9 & (     11,      1)                   &      10 &      10 &      10 &      10 &      10 & \textbf{(     10,      0)} &      10 &      10 &      10 &      10 &      10 & (      9,      1)          &      10 &      10 &      10 &      10 &      10 & (      5,      4)          &      10 &      10 &      11 &      11 &      11\\
(      6,      6)          &      11 &      11 &      10 &      10 &      10 & \textbf{(      9,      7)}          &      11 &      11 &      11 &      11 &      11 & \textbf{(      8,      6)} &      11 &      11 &      11 &      11 &      11 & (      9,      3)          &      11 &      11 &      11 &      11 &      11 & (      6,      3)          &      11 &      11 &      10 &      10 &      10\\
(      8,      4)          &      13 &      13 &      13 &      13 &      13 & \textbf{(     11,      3)}          &      12 &      12 &      12 &      12 &      12 & (     10,      2)          &      12 &      12 &      12 &      13 &      13 & (      7,      7)          &      12 &      12 &      12 &      12 &      12 & (      7,      0)          &      12 &      12 &      12 &      12 &      12\\
\textbf{(      ?,      ?)} &      13 &      14 &      14 &      14 &      14 & (     11,      5)                   &      13 &      13 &      13 &      13 &      13 & (     10,      4)          &      13 &      13 &      14 &      14 &      14 & (      9,      5)          &      13 &      13 &      13 &      13 &      13 & (      7,      2)          &      13 &      13 &      13 &      13 &      13\\
\textbf{(      ?,      ?)} &      14 &      13 &      13 &      13 &      13 & \textbf{(     13,      1)}          &      14 &      15 &      15 &      15 &      15 & (     10,      6)          &      14 &      14 &      13 &      12 &      12 & (     11,      1)          &      14 &      14 &      14 &      14 &      14 & (      6,      5)          &      14 &      14 &      14 &      14 &      14\\
(     10,      2)          &      15 &      15 &      15 &      15 &      15 & \textbf{(     11,      7)}          &      15 &      14 &      14 &      14 &      14 & (     12,      0)          &      15 &      15 &      15 &      15 &      15 & \textbf{(     11,      3)} &      15 &      15 &      15 &      15 &      15 & \textbf{(      7,      4)} &      15 &      15 &      15 &      15 &      15\\
(     10,      4)          &      16 &      16 &      18 &      19 &      19 & (     13,      3)                   &      16 &      16 &      16 &      16 &      16 & (     12,      2)          &      16 &      16 &      16 &      16 &      16 & \textbf{(      9,      7)} &      16 &      16 &      16 &      16 &      16 & \textbf{(      8,      1)} &      16 &      16 &      16 &      16 &      16\\
(      8,      8)          &      17 &      17 &      16 &      16 &      16 & (     13,      5)                   &      17 &      17 &      17 &      17 &      17 & (     12,      4)          &      17 &      17 &      18 &      18 &      18 & (     11,      5)          &      17 &      17 &      17 &      17 &      17 & (      8,      3)          &      17 &      17 &      17 &      17 &      17\\
(     10,      6)          &      18 &      18 &      17 &      17 &      17 & (     11,      9)                   &      18 &      18 &      18 &      19 &      19 & (     10,      8)          &      18 &      18 &      17 &      17 &      17 & (      9,      9)          &      18 &      18 &      18 &      18 &      18 & (      9,      0)          &      18 &      18 &      18 &      18 &      18\\
(     12,      0)          &      19 &      19 &      19 &      18 &      18 & (     13,      7)                   &      19 &      19 &      20 &      20 &      20 & (     12,      6)          &      19 &      19 &      19 &      19 &      19 & \textbf{(     11,      7)} &      19 &      20 &      20 &      20 &      20 & \textbf{(      7,      6)} &      19 &      19 &      19 &      19 &      19\\
(     12,      2)          &      20 &      20 &      26 &      26 &      26 & (     15,      1)                   &      20 &      20 &      19 &      18 &      18 & (     14,      0)          &      20 &      20 &      20 &      20 &      20 & \textbf{(     13,      1)} &      20 &      19 &      19 &      19 &      19 & \textbf{(      9,      2)} &      20 &      21 &      21 &      21 &      21\\
(     12,      4)          &      21 &      21 &      20 &      20 &      20 & (     15,      3)                   &      21 &      21 &      21 &      21 &      21 & (     14,      2)          &      21 &      21 &      21 &      21 &      21 & (     13,      3)          &      21 &      21 &      21 &      21 &      21 & (      8,      5)          &      21 &      20 &      20 &      20 &      20\\
(     10,      8)          &      22 &      22 &      21 &      21 &      21 & \textbf{(     13,      9)}          &      22 &      22 &      22 &      22 &      22 & (     12,      8)          &      22 &      22 &      22 &      22 &      22 & (     13,      5)          &      22 &      22 &      22 &      22 &      22 & (      9,      4)          &      22 &      22 &      22 &      22 &      22\\
(     12,      6)          &      23 &      23 &      22 &      22 &      22 & \textbf{(     15,      5)}          &      23 &      23 &      23 &      23 &      23 & (     14,      4)          &      23 &      23 &      23 &      23 &      23 & (     11,      9)          &      23 &      23 &      24 &      23 &      23 & (     10,      1)          &      23 &      23 &      23 &      23 &      23\\
(     14,      0)          &      24 &      24 &      23 &      23 &      23 & (     15,      7)                   &      24 &      24 &      26 &      26 &      26 & (     14,      6)          &      24 &      24 &      24 &      24 &      24 & (     13,      7)          &      24 &      24 &      23 &      26 &      26 & (     10,      3)          &      24 &      24 &      26 &      26 &      26\\
                           &         &         &         &         &         & \textbf{(     17,      1)}          &      25 &      25 &      25 &      25 &      25 & (     12,     10)          &      25 &      25 &      25 &      25 &      25 & (     15,      1)          &      25 &      25 &      25 &      24 &      24 & (      8,      7)          &      25 &      25 &      24 &      24 &      24\\
                           &         &         &         &         &         & \textbf{(     13,     11)}          &      26 &      26 &      27 &      27 &      27 &                            &         &         &         &         &         & (     15,      3)          &      26 &      26 &      26 &      25 &      25 &                            &         &         &         &         &        \\
                           &         &         &         &         &         & (     17,      3)                   &      27 &      27 &      24 &      24 &      24 &                            &         &         &         &         &         & (     11,     11)          &      27 &      27 &      27 &      27 &      27 &                            &         &         &         &         &        \\
                           &         &         &         &         &         & (     15,      9)                   &      28 &      28 &      28 &         &         &                            &         &         &         &         &         &                            &         &         &         &         &         &                            &         &         &         &         &        \\

\end{tabular}
}
\caption{
Table of Eigenfunction Correspondences:
For the Square, $(a,b)\leftrightarrow \cos(\pi ax)\cos(\pi by)\pm\cos(\pi ay)\cos(\pi bx)$,
($+$ in columns $1,4$, $-$ in columns $2,3$).
Two-dimensional eigenspaces in bold.
Question marks denote blatantly non-canonical splitting of a two-dimensional eigenspace
}
\label{tabletwo}
\end{table}

\clpg

\clpg

\clpg

\clpg

    \clpg

    \clpg

    \clpg

    \begin{figure}[htbp!]
    \begin{center}
    \includegraphics[scale=1]{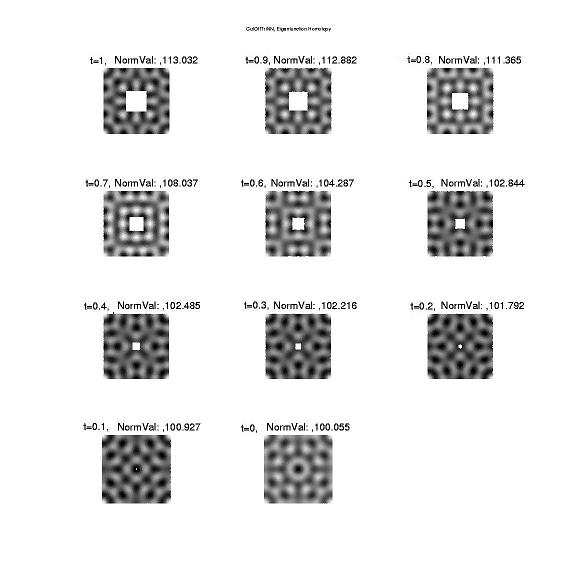}
    \end{center}
    \caption{$(1++)$ EigVal 13}
    \label{figtwentyfive}
    \end{figure}

    \begin{figure}[htbp!]
    \begin{center}
    \includegraphics[scale=1]{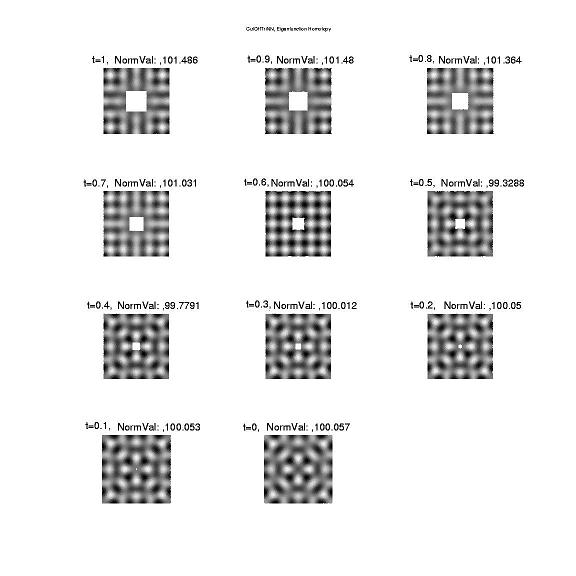}
    \end{center}
    \caption{$(1++)$ EigVal 14}
    \label{figtwentyfive.1}
    \end{figure}

\clearpage


\bibliographystyle{plain}

\end{document}